\documentclass[final,3p]{elsarticle}

\usepackage{amssymb}
\usepackage{amsmath}
 \usepackage{amsthm}
 \usepackage{subfig}
 \usepackage{gensymb}

\usepackage{lineno}
\usepackage[colorlinks=true,breaklinks=true,pdftex]{hyperref}
\modulolinenumbers[1]

\biboptions{authoryear}

\begin{document}

\begin{frontmatter}

\title{Adaptive optimization of isogeometric multi-patch discretizations using artificial neural networks}

\author[jku,ricam]{Dany Rios}
\author[jku]{Felix Scholz}
\ead{felix.scholz@jku.at}
\author[ricam]{Thomas Takacs}
\ead{thomas.takacs@ricam.oeaw.ac.at}

\address[jku]{Johannes Kepler University Linz, Austria}
\address[ricam]{RICAM, Austrian Academy of Sciences}

\date{March 25, 2024}

\begin{abstract}
In isogeometric analysis, isogeometric function spaces are employed for accurately representing the solution to a partial differential equation (PDE) on a parameterized domain.
They are generated from a tensor-product spline space  by composing the basis functions with the inverse of the parameterization.
Depending on the geometry of the domain and on the data of the PDE, the solution might not have maximum Sobolev regularity, leading to a reduced convergence rate.
In this case it is necessary to reduce the local mesh size close to the singularities.
The classical approach is to perform adaptive $h$-refinement, which  either leads  to an unnecessarily large number of degrees of freedom or to a spline space that does not possess a tensor-product structure.
Based on the concept of $r$-adaptivity we present a novel approach for finding a suitable isogeometric function space for a given PDE without sacrificing the tensor-product structure of the underlying spline space.
In particular, we use the fact that different reparameterizations of the same computational domain lead to different isogeometric function spaces while preserving the geometry.
Starting from a multi-patch domain consisting of bilinearly parameterized patches, we aim to find the biquadratic multi-patch parameterization that leads to the isogeometric function space with the smallest best approximation error of the solution.
In order to estimate the location of the optimal control points, we employ a trained residual neural network that is applied to the graph surfaces of the approximated solution and its derivatives.
In our experimental results, we observe that our new method results in a vast improvement of the approximation error for different PDE problems on multi-patch domains.
\end{abstract}


\end{frontmatter}

\section{Introduction}

In this paper we propose a new approach to adaptively optimize an isogeometric discretization. Isogeometric Analysis (IGA), see~\cite{hughes2005isogeometric,beiraodaveiga2014mathematical}, is a numerical method that uses NURBS or B-spline representations both for the geometry and for the discretization space. We focus on planar domains that are segmented into quadrilateral B\'ezier patches. The basic idea is to adapt the parameterization of the domain in order to reduce the approximation order of the isogeometric discretization. As a model problem we consider the Poisson equation with Dirichlet boundary conditions, however, the approach may be applied directly to more general problems, since it does not rely on the underlying PDE. In the following we summarize existing methods of adaptivity in IGA and then give an overview of our approach and motivate its set-up.

For problems that have singular solutions or solutions of low Sobolev regularity, the solution cannot be resolved well near the singularity and the discretization error is large there. Thus it is better to adaptively refine near the singularity instead of performing global refinement. In IGA, many approaches for $h$-refinement exist, which are based on locally refinable splines, such as T-splines, LR B-splines or THB-splines (see~\cite{bazilevs2010isogeometric,dorfel2010adaptive},~\cite{johannessen2014isogeometric} and~\cite{giannelli2016thb}, respectively). See also the comparison paper~\cite{hennig2017adaptive}. A neural-network based approach for finding a locally refined mesh was proposed in \cite{chan2022locally}. To reduce the error, one may also increase the degree locally, which is called $p$-refinement, or $k$-refinement, if the spline regularity is raised simultaneously, cf.~\cite{cottrell2007studies}. It is also possible to combine $p$- and $h$-refinement, as e.g. in~\cite{schillinger2011unfitted,liu2016hybrid,kamber20222}. 
When using NURBS-based parameterizations, one can also update the NURBS-weights of the discretization space to improve the discretization, as in~\cite{taheri2020adaptive}.

Alternatively, one may also reduce discretization errors by adapting the parameterization of the domain. This approach was termed $r$-refinement in~\cite{xu2011parameterization} and investigated further in several studies~\cite{xu2019efficient,ji2023curvature}. A similar approach based on optimal transport was developed in~\cite{bahari2024adaptive}. One key advantage of this type of refinement is that the spaces remain tensor-product, which implies a much simpler data structure in contrast to locally refined splines and allows for very efficient numerical quadrature and assembly, cf.~\cite{scholz2018partial,bressan2019sum}. However, existing $r$-refinement approaches are costly, since they require solving an optimization problem to obtain a suitable parameterization. Instead, we want to find a good update of the parameterization using a neural network. Such an approach is cheap and can, in principle, yield near-optimal results for many applications. We moreover explore the possibility of combining such an $r$-refinement approach with standard, global $h$-refinement. A similar study was performed in~\cite{basappa2016adaptive}.

The basic set-up of our approach is as follows. We assume that the domain is split into quadrilaterals. First we solve the PDE problem that we are interested in on the initial geometry. We then apply an artificial neural network to find a biquadratic reparameterization of each patch, using as input the initial solution of the PDE (evaluated at sampled points). This reparameterization should then yield a smaller discretization error than the initial parameterization.

The paper is organized as follows. In Section~\ref{sec:isogeometric-spaces} we introduce the isogeometric discretization spaces that we consider throughout the paper. In Section~\ref{sec:modelProblems} we present the two model problems that are considered in this work. In Section~\ref{sec:overview-method} we give an overview of the proposed method, which is then tested on several examples in Section~\ref{sec:numerical-experiments}. Possible extensions and conclusions are discussed in Section~\ref{sec:extension-conclusions}.

\section{Isogeometric discretization spaces}\label{sec:isogeometric-spaces}
We consider planar domains $\Omega\subset\mathbb{R}^2$, which are segmented into quadrilateral patches $\Omega^k$, i.e.,
\[
 \overline\Omega = \bigcup_{k=1}^K \overline{\Omega^k},
\]
where $\partial \Omega^k \cap \partial\Omega^{k'}$, for $k\neq k'$, is empty, a vertex or an entire edge of both quadrilaterals. Each subdomain is parameterized by 
\[
 G^k: \widehat\Omega \rightarrow \Omega^k,
\]
where $\widehat\Omega = \left]0,1\right[^2$ and $G^k$ is a regular mapping ($\det\nabla G^k \geq \underline{c} >0$). For simplicity we use the notation $G=(G^1,G^2,\ldots,G^K)$.

Let $\mathbb{P}^1_p$ be the space of univariate polynomials of degree $p$ and let $\mathbb{P}^2_p$ and $\mathbb{Q}^2_p$ be the spaces of bivariate polynomials of total degree $p$ and maximum degree $p$, respectively. Let moreover $S^d_{p,h}$ be the $d$-variate B-spline space of degree $p$, regularity $p-1$ and uniform mesh size $h$ in each direction. We have $S^2_{p,h} = S^1_{p,h} \otimes S^1_{p,h}$.

For simplicity, we assume that the initial patch parameterizations $G^k_{init}$ are bilinear, i.e., $G^k_{init} \in \mathbb{Q}^2_1 \times \mathbb{Q}^2_1$. We are interested in biquadratic reparameterizations $G^k \in \mathbb{Q}^2_2 \times \mathbb{Q}^2_2$, which minimize the approximation error as explained later.

\begin{figure}[ht]
\begin{center}
 	\subfloat[Uniform refinement of the isogeometric function space]{\includegraphics[width=.21\textwidth, viewport=10cm 3cm 31cm 23cm, clip]{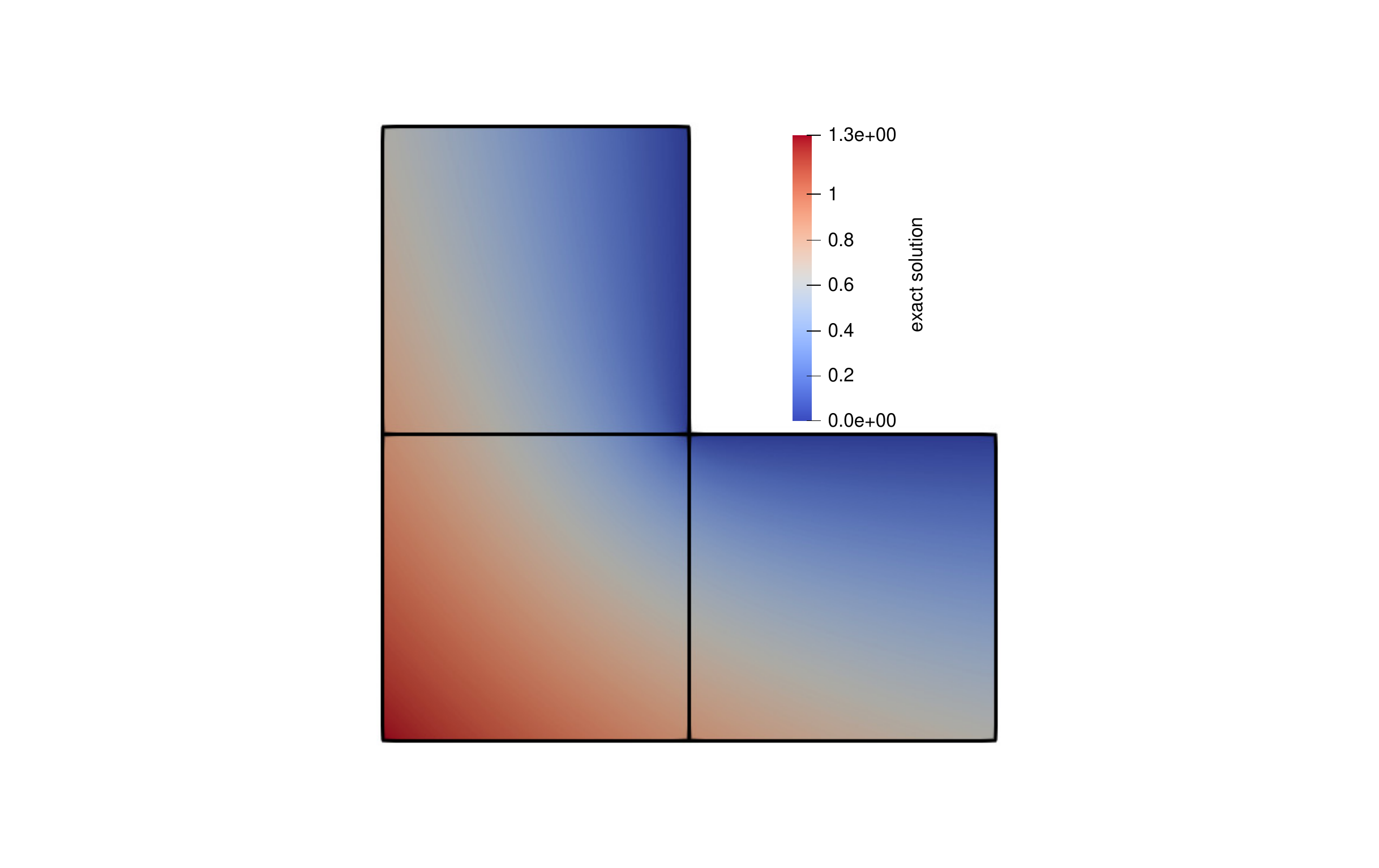}
	\includegraphics[width=.21\textwidth, viewport=10cm 3cm 31cm 23cm, clip]{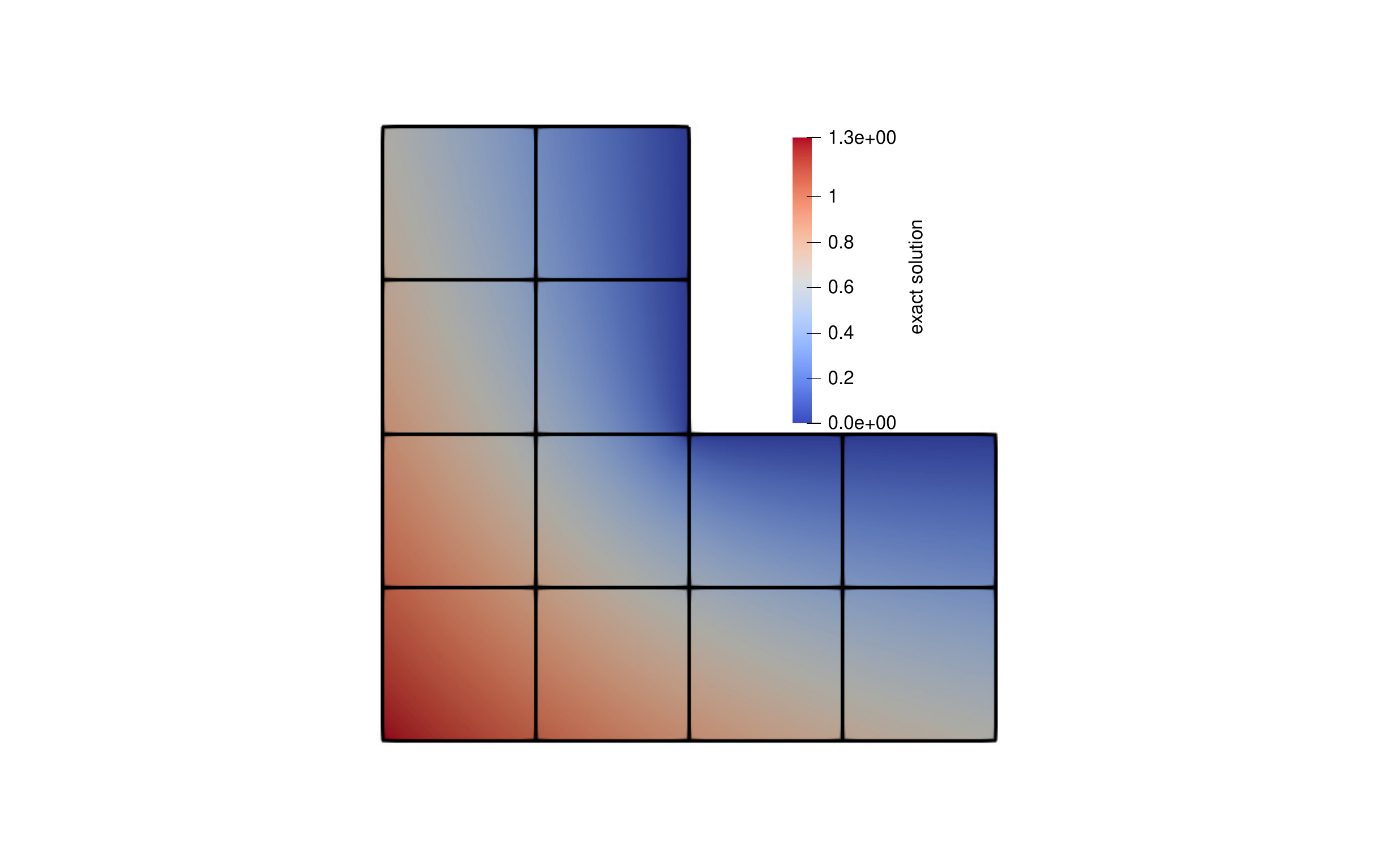}
 	\includegraphics[width=.21\textwidth, viewport=10cm 3cm 31cm 23cm, clip]{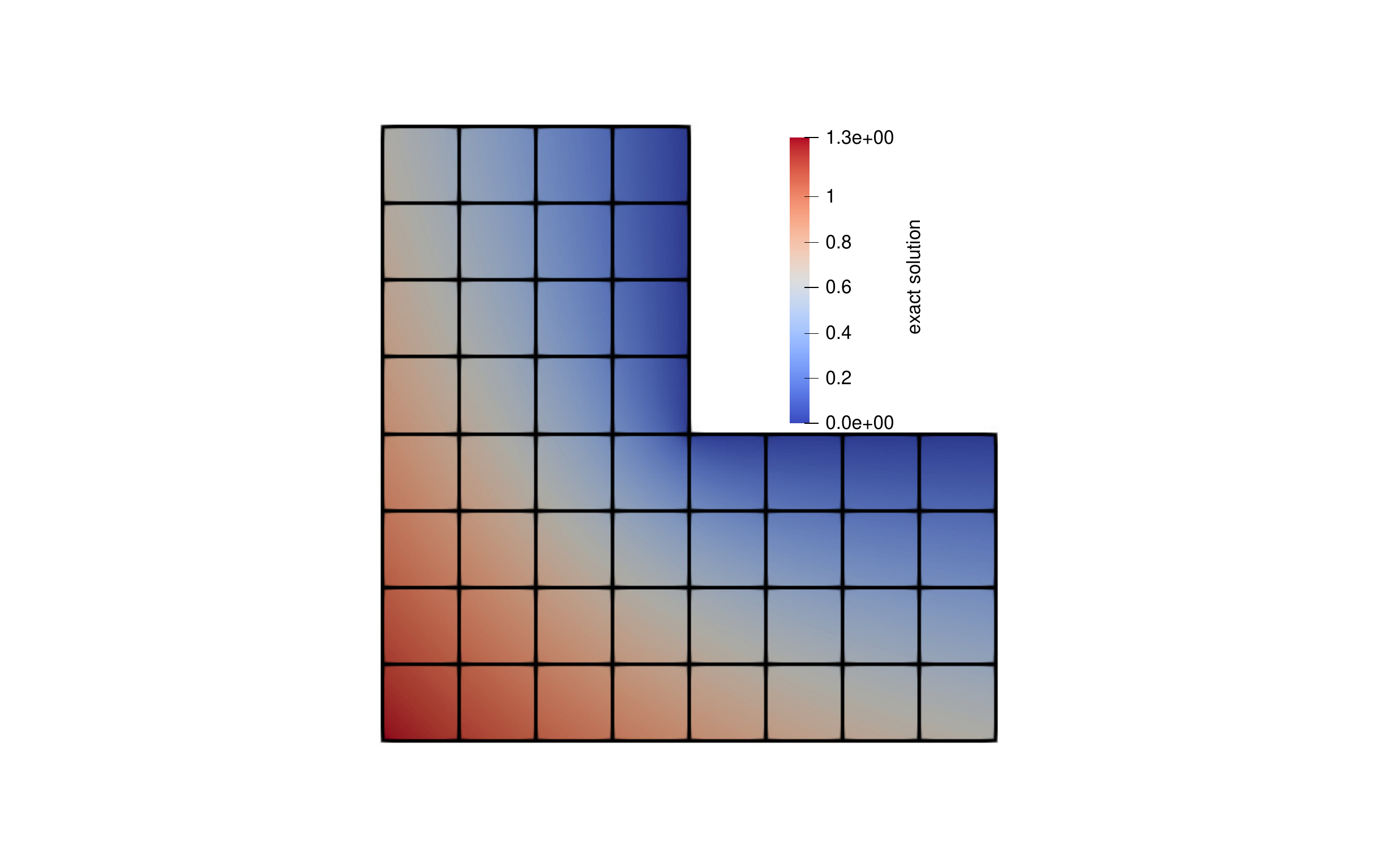}} \\
 	\subfloat[Optimized refinement of the isogeometric function space]{\includegraphics[width=.21\textwidth, viewport=10cm 3cm 31cm 23cm, clip]{exactLshapeStationaryHeatFlat-eps-converted-to}
 \includegraphics[width=0.21\textwidth, viewport=10cm 3cm 31cm 23cm, clip]{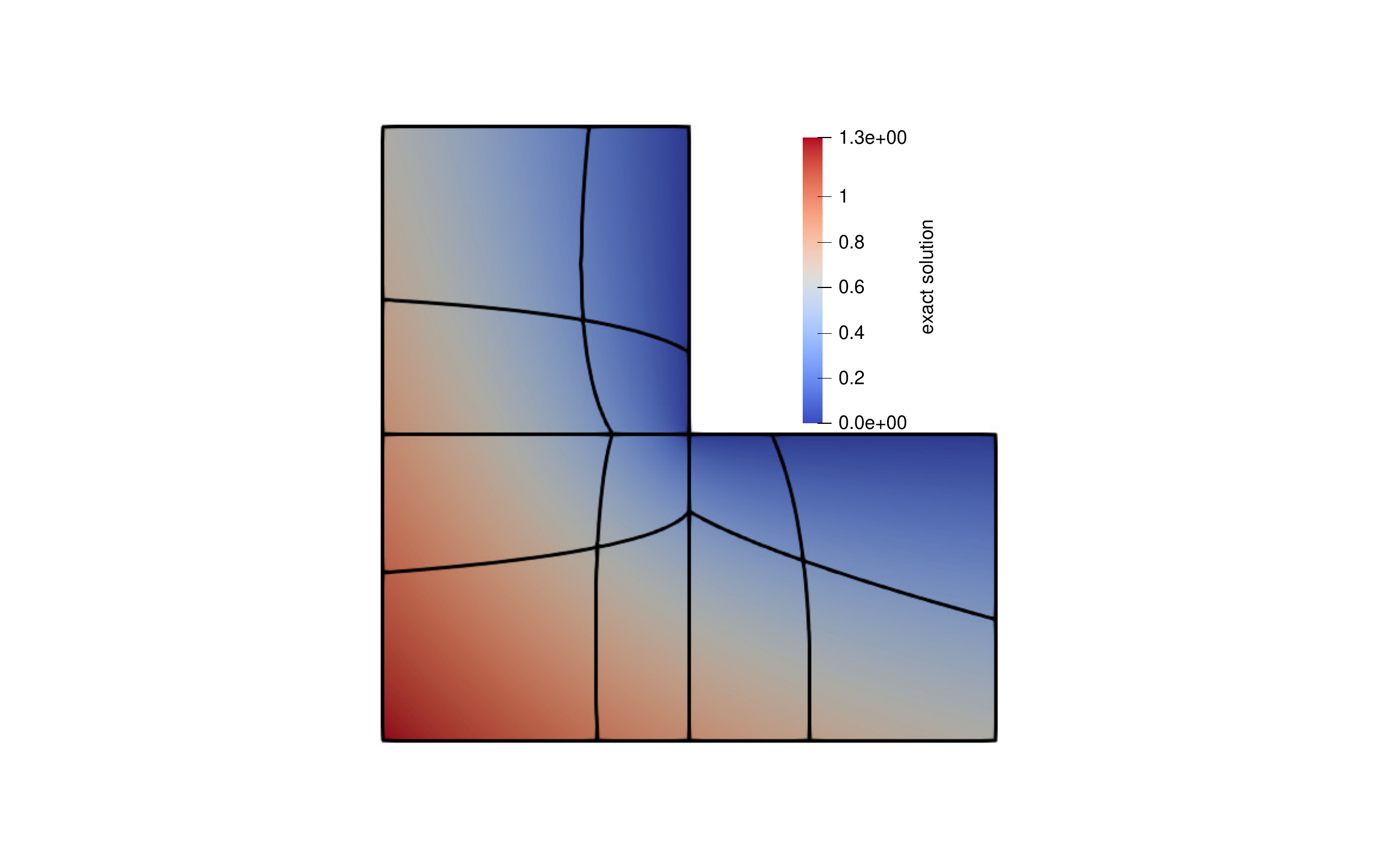}
 \includegraphics[width=0.21\textwidth, viewport=10cm 3cm 31cm 23cm, clip]{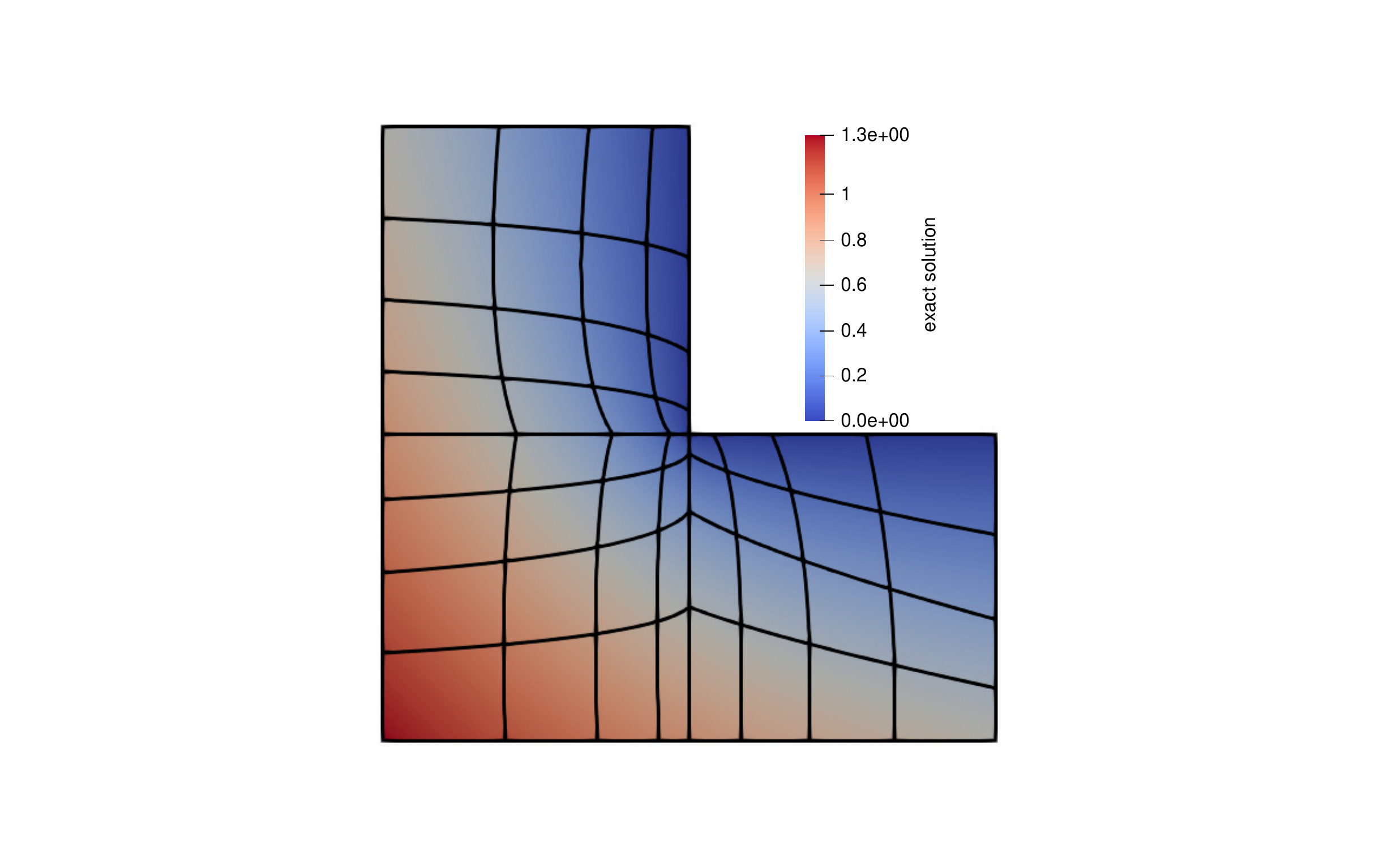}}
 \end{center}
 \caption{Stationary solution to the heat equation on an L-shaped domain}\label{fig:heat-Lshape}
\end{figure}

The isogeometric function spaces are defined as 
\[
 V_{h,G} = \{ \varphi \in C^0(\Omega) : \varphi \circ G^k \in S^2_{p,h}, \forall k \}.
\]
We observe that different parameterizations $G$ of the same domain lead to different isogeometric function spaces. See Figure~\ref{fig:heat-Lshape} for a comparison of different isogeometric discretizations. The different convergence rates are shown in Figure~\ref{fig:L-shape-heat-all}.

\section{Model problems}\label{sec:modelProblems}

In this paper we consider the Poisson equation as a model problem. It is formulated as finding $u$ such that
\[
\begin{array}{rcll}
 -\Delta u &=& f &\quad\mbox{ in }\Omega \\
 u &=& g_D &\quad\mbox{ on }\Gamma_D \\
 n \cdot \nabla u &=& g_N &\quad\mbox{ on }\Gamma_N,
\end{array}
\]
with $\partial\Omega = \Gamma_D \cup \Gamma_N$, $n$ being the outward pointing normal vector to $\Gamma_N$ and $f$, $g_D$ and $g_N$ are suitably regular data of the PDE. This problem has the variational form: Find $u\in H^1_{g_D} (\Omega) = \{ u \in H^1(\Omega) : u = g_D \; \mbox{ on } \Gamma_D\}$ with
\begin{linenomath}
\begin{equation}\label{eq:poisson-variational}
 a(u,v)= F(v) \quad\mbox{ for all }v\in H^1_{0} (\Omega),
\end{equation}
\end{linenomath}
where
\[
 a(u,v) := \int_\Omega \nabla u \nabla v \mathrm{d}x 
\]
and
\[
 F(v) := \int_\Omega f v \mathrm{d}x + \int_{\Gamma_N} g_N v \mathrm{d}x.
\]
The discrete problem that we solve is then given by finding $u_h\in V_{h,G,g_D} = V_{h,G} \cap H^1_{g_D} (\Omega)$ with
\begin{linenomath}\begin{equation}\label{eq:poisson-discrete}
 a(u_h,v_h)= F(v_h) \quad\mbox{ for all }v_h\in V_{h,G,0}.
\end{equation}\end{linenomath}
Moreover, when initially testing our approach we solve an $L^2$-fitting problem 
\[
 \min_{u_h \in V_{h,G}}\| u_h - u \|^2_{L^2(\Omega)},
\]
which is given in variational form as 
\[
 \int_\Omega u_h v_h \mathrm{d}x = \int_\Omega u v_h \mathrm{d}x \quad\mbox{ for all }v_h\in V_{h,G}.
\]

For both the Poisson problem and the $L^2$-fitting problem we assume that the exact solution $u$ does not possess a high order of Sobolev regularity, i.e., that $u\in H^1(\Omega)$ or $u\in L^2(\Omega)$, respectively, but $u\notin H^{p+1}(\Omega)$. Depending on the Sobolev regularity of the exact solution, a suitable mesh grading is needed.

We are now interested in the discretization error
\[
 | u - u_h |_{H^1(\Omega)} = \| \nabla u - \nabla u_h \|_{L^2(\Omega)}
\]
or
\[
 \| u - u_h \|_{L^2(\Omega)}.
\]
We intend to minimize the error by changing the parameterization $G$ of the domain and thus changing the discretization space $V_{h,G}$, which depends non-linearly on $G$.

\section{Description of our method}\label{sec:overview-method}

As described in the previous section, we assume to have given a domain $\Omega$ and a model problem~\eqref{eq:poisson-discrete} over the domain, which has the solution $u_h \in V_{h,G,g_d}= V_{h,G} \cap H^1_{g_D} (\Omega)$. It is our goal to minimize the error with respect to the parameterization $G$, i.e.,
\[
 \| u - u_h \|_V \rightarrow \min_{G},
\]
where $u \in V$ is the solution of~\eqref{eq:poisson-variational}.

This minimization problem is difficult as the discretization spaces $V_{h,G}$ and discrete solutions $u_h$ depend non-linearly on $G$. Thus we propose an approach based on machine learning to optimize for $G$.
\subsection{General overview of the method}\label{general}
To be able to apply a neural network to general domain compositions, we need to segment the problem. We consider the following workflow: We first split the domain into quadrilaterals and find a biquadratic reparameterization for each quadrilateral. Then we use an existing, already trained, artificial neural network, which was developed by~\cite{pointcloudParam}, to reparameterize each quadrilateral, making sure that the reparameterizations match at patch interfaces.

The network is defined on triangular patches. It takes as input (unsorted) point data and finds parameter values to the given data points. It was trained in such a way that it optimizes the parameter values such that the error of the $\ell^2$-projection to the data points is minimized. For most configurations, the network predicts a good approximation to the original parameter values if the point data was sampled from a quadratic surface, cf.~\cite{pointcloudParam}. The way the network is set up, it can be used to determine the best approximating $\ell^2$-fit of a quadratic surface to the given point data.

The use of such a neural network has several advantages. Considering only point data, the network can be applied to a number of different discretizations and is not restricted to a specific set-up. The network yields parameter values with respect to a triangular patch, since the corresponding parameterization problem on a triangle is well-defined. A best approximating triangular surface can be used as a good initial guess for optimizing a tensor-product polynomial or spline parameterization. This is based upon two assumptions that we verify in our experiments.

\paragraph{Assumption 1} When splitting a quadrilateral into triangles, the averaged reparameterizations of the best approximating quadratic surfaces yield a biquadratic surface which is close to the best approximating biquadratic surface.

\bigskip

Our numerical evidence suggests that we can use properly averaged updates for triangular patches derived from the network.
To optimize the results, we tested and compared several different averaging strategies.

\paragraph{Assumption 2} Taking the best approximating surface to sampled point data and projecting this surface onto the $(x,y)$-plane, i.e., onto the domain $\Omega$, we obtain a reparameterization that is close to the optimal reparameterization, which yields the smallest discretization error.

\bigskip

In a general notation, we want to minimize the projection error in the energy norm:
\[
 \| ({\rm I}-\Pi_{h}^G) (u) \|_{V} \rightarrow \min_{G},
\]
where $u$ is the unknown solution, $\Pi_{h}^G$ is a projector into $V_{h,G}$, and we minimize with respect to the parameterization $G$, such that $G(\widehat{\Omega}) = \Omega$ and $G(\partial\widehat{\Omega}) = \partial\Omega$.
In case of the $L^2$-fitting problem,  $\|.\|_V$ is the $L^2$-norm and $\Pi_{h}^G$
is the $L^2$-orthogonal projector.
In case of the Poisson problem, $\|.\|_V$ is the $H^1$-seminorm and $\Pi_{h}^G$ is the $H^1$-orthogonal projector.

Since the network is based on an $L^2$-approximation, we always approximate the $L^2$-projection and we optimize the $L^2$-error, i.e., we consider the problem
\begin{linenomath}\begin{equation}\label{eq:L2-error}
 \| ({\rm I}-\Pi_{h}^G) (u) \|_{L^2(\Omega)} \rightarrow \min_{G}.
\end{equation}\end{linenomath}
However, in case of the Poisson problem, this can be replaced by an approximation of the partial derivatives of the unknown solution $u$, which should mimic an $H^1$-projection. This corresponds to solving the problem
\begin{linenomath}\begin{equation}\label{eq:L2-derivative-error}
 \left\| ({\rm I}-\Pi_{h}^G) \left(\frac{\partial u}{\partial x}\right) \right\|^2_{L^2(\Omega)}+ \left\| ({\rm I}-\Pi_{h}^G) \left(\frac{\partial u}{\partial y}\right) \right\|^2_{L^2(\Omega)} \rightarrow \min_{G}.
\end{equation}\end{linenomath}

Moreover, instead of minimizing a norm directly, which would require complete knowledge of the function $u$, we sample points of an approximation of $u$, which are then fed to the network. In the numerical experiments we compare different ways to sample points, depending on the underlying model problem. This is discussed in more detail in Section~\ref{sec:point-sampling}. 

\subsection{Description of the employed neural network}\label{NN}
We use the point cloud parameterization network from \cite{pointcloudParam}. It is a residual neural network with fixed input and output dimensions equal to 36. The input consists of twelve standardized points in $\mathbb R^3$ and its output is interpreted as barycentric coordinates for each of the input points, with respect to the standard triangle with vertices $(0,0)^T$, $(1,0)^T$ and $(\frac12,\frac{\sqrt 3}{2})^T$.

The architecture consists of a linear input layer, four residual blocks and a linear output layer. The activation function of all hidden layers is ReLU and only the output activation function is the sigmoid function in order to ensure that the output values lie between zero and one. Finally, the output barycentric coordinates for each input point are enforced to sum up to 1.

As described in \cite{pointcloudParam}, the network was trained on a dataset consisting of 500,000 point clouds, sampled from quadratic B\'ezier surfaces in an unsupervised way, with the loss function representing the $\ell^2$-error when fitting the point cloud with a quadratic surface using the the barycentric coordinates that were predicted by the network.
Before applying the network to a point cloud of twelve points that is to be parameterized with respect to an arbitrary triangle, the point cloud is standardized by an affine transformation that maps the triangle's vertices to the standard triangle.

In \cite{pointcloudParam}, the output of the neural network is used as an initialization for different non-linear optimization methods, of which the Levenberg--Marquardt algorithm performs the best. However, in the context of our present problem, we observed that using the network's predictions directly without further optimization leads to a much more robust method. The reason is that the Levenberg--Marquardt algorithm often converges to parameterizations with parameters that lie outside of the triangular parameter domain, i.e., barycentric coordinates with negative components.
While these parameterizations are minima of the fitting problem, they are not suitable for optimizing isogeometric function spaces using our method since our solution is only defined inside the parameter domain. Because of the network architecture,
the network consistently outputs non-negative barycentric coordinates.

\subsection{Processing the network output}
While the employed network was originally designed to parameterize data over triangular domains, we apply it in this work to our problem of finding optimal reparameterizations of quadrilateral multi-patch domains. We do this by splitting the problem into triangles, finding optimized triangular B\'ezier reparameterizations and performing a suitable averaging procedure, based on the two assumptions made in Section~\ref{general}.
\subsubsection{Optimizing quadratic parameterizations of triangular domains}\label{sec:triangle}
We first describe how to apply the parameterization network to optimize the parameterization of a triangle. Let $\Delta$ be an arbitrary triangle and let $u:\Delta\rightarrow\mathbb R$ be the function to be approximated in an isogeometric (isoparametric) function space over $\Delta$.

As an initial parameterization, we choose the linear parameterization represented as a quadratic B\'ezier patch
\begin{linenomath}$$T_{\mathit{init}}(\alpha, \beta, \gamma) = \sum_{\substack{i+j+k=2 \\ i,j,k\geq 0}} T_{ijk}^{\mathit{init}} B_{ijk}^2(\alpha, \beta, \gamma),$$\end{linenomath}
where $\alpha,\beta,\gamma\geq 0$ are barycentric coordinates with $\alpha+\beta+\gamma=1$, $B_{ijk}^2$ are the quadratic Bernstein polynomials, the control points $T_{200}^{\mathit{init}}, T_{020}^{\mathit{init}}, T_{002}^{\mathit{init}}\in\mathbb R^2$ are the vertices of $\Delta$ and the edge control points $T_{110}^{\mathit{init}}, T_{101}^{\mathit{init}}, T_{011}^{\mathit{init}}$ correspond to the midpoints of the edges of $\Delta$.

In order to optimize the control points of $T$, the neural network is applied to sets of 12 points sampled from the graph surface.
More precisely, we create a number of $M$ sets $Q_r$, each consisting of 12 points
\begin{linenomath}\begin{equation*}
 P^r_i = \left(T_{\mathit{init}}(\alpha^r_i, \beta^r_i, \gamma^r_i), u\circ T_{\mathit{init}}(\alpha^r_i,\beta^r_i,\gamma^r_i)\right).
\end{equation*}\end{linenomath}
generated from randomly sampled barycentric coordinates $\alpha^r_i, \beta^r_i, \gamma^r_i\geq0$ with $\alpha^r_i+\beta^r_i+\gamma_i^r = 1$.

After standardizing each set $Q_r$ as described in Section~\ref{NN}, the neural network is applied in order to predict optimal barycentric coordinates $\tilde\alpha_i^r,\tilde\beta_i^r, \tilde\gamma_i^r$.

In order to realize each parameterization of $\Delta$ that was predicted by the neural network, we compute a reparameterization $T^r$ by solving the linear least squares problem
\[
 \min_{\delta_{011}, \delta_{101}, \delta_{110}\in[0,1]}\sum_{i=1}^{12} \|T^r(\tilde\alpha_i^r, \tilde\beta_i^r, \tilde\gamma_i^r) - T(\alpha^r_i, \beta_i^r,\gamma_i^r)\|^2,
\]
where $T^r$ is the quadratic B\'ezier parameterization with vertex control points
\begin{linenomath}\begin{equation*}
 T^r_{002} = T^{\mathit{init}}_{002},\quad T^r_{020} = T^{\mathit{init}}_{020},\quad  T^r_{200} = T^{\mathit{init}}_{200}
\end{equation*}\end{linenomath}
and with edge control points
\begin{linenomath}\begin{equation*}
 T^r_{011} = (1-\delta_{011})T^r_{002} + \delta_{011}T^r_{020},\quad  T^r_{101} = (1-\delta_{101})T^r_{002} + \delta_{101}T^r_{200},\quad  T^r_{110} = (1-\delta_{110})T^r_{020} + \delta_{110}T^r_{200}.
\end{equation*}\end{linenomath}
This enforces that $T^r$ is the quadratic parameterization of $\Delta$ that best approximates the parameterization $\tilde\alpha_i^r,\tilde\beta_i^r, \tilde\gamma_i^r$ predicted by the neural network. From  the reparameterizations $T^r$ for $r=1,\ldots, M$ we obtain a single reparameterization $T$ by averaging the control points.

We iterate this process by sampling from the graph surface of $u$ using the deformed parameterization $T$ until a suitable stopping criterion is satisfied. 
Using batches of size $M>1$ in the optimization process leads to an improved stability compared to applying the method to a single sampled point cloud.

\subsubsection{Optimizing biquadratic parameterizations of quadrilateral domains}\label{biquadratic}
We now extend this method to quadrilateral domains.
To this end, we assume to have given a bilinear quadrilateral $Q$, represented by a biquadratic B\'ezier patch
\begin{linenomath}\begin{equation*}
 G_{init}(s,t) = \sum_{i=0}^2 \sum_{j=0}^2 G^{init}_{ij}\hat B_{ij}^{(2,2)}(s,t),
\end{equation*}\end{linenomath}
where $\hat B_{ij}^{(2,2)}$ are the biquadratic tensor-product Bernstein polynomials and the control points $G^{init}_{ij}$ are computed from the corners $G_{00}$, $G_{02}$, $G_{20}$, $G_{22}$ of the quadrilateral by degree elevation. We want to find a biquadratic parameterization $G$ of the same quadrilateral $Q$ that minimizes the approximation error to the target function $u: Q\rightarrow\mathbb R$. Thus the corner control points of $G$ are equal to the corners of $Q$ and the boundary segments remain straight line segments, parameterized non-uniformly, e.g., $G_{10} \in \overline{G_{00}G_{20}}$.

\begin{figure}[h!]
	\centering
	\includegraphics{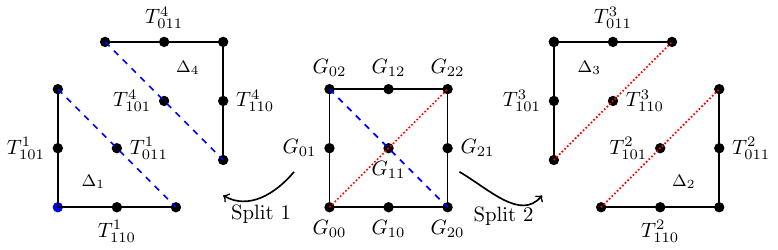}
	\caption{Splitting a quadrilateral domain into four triangles.}
	\label{fig:Quad_Splitting}
\end{figure}

To find such a reparameterization, we use the method for triangular domains described in the previous section. We split the quadrilateral $Q$ into four triangles $\Delta_1 = G_{00}G_{20}G_{02}$, $\Delta_2 = G_{20}G_{22}G_{00}$, $\Delta_3 = G_{02}G_{00}G_{22}$ and $\Delta_4 = G_{22}G_{02}G_{20}$. To each of these triangles we apply the previously described algorithm to obtain optimized quadratic reparameterizations $T^w$ of $\Delta_w$, with control points $T^w_{110}$, $T^w_{101}$ and $T^w_{011}$, for $w=1,2,3,4$, see Fig.~\ref{fig:Quad_Splitting}.

In order to compute the biquadratic parameterization $G$ for the quadrilateral $Q$ from the optimized reparameterizations $T^w$ of the triangles $\Delta_w$, we compare three different averaging procedures.

\paragraph{Best-approximating inverse}

If we assume to have given a quadrilateral mapping $G$, one can directly compute the best approximating triangular parameterizations $T^w$. In the first approach we compute the unknown biquadratic deformation by looking at the pseudo-inverse of this approximation procedure.

Assuming to have given $G$, we approximate it with quadratic mappings on all four triangles, i.e., we find the approximations $T^{w,\mathit{appr}}$ by $L^2$-fitting
\[
 \| G|_{\Delta_w} - T^{w,\mathit{appr}} \|_{L^2(\Delta_w)} \rightarrow \min.
\]
This yields a linear mapping $\mathcal{A}$ on the control points
\[
 \mathcal{A}\left( (G_{ij})_{ij=0}^2 \right) = \left( ((T^{w,\mathit{appr}}_{ijk})_{i+j+k=2})_{w=1}^4 \right).
\]

Of this linear operator we compute the Moore--Penrose pseudoinverse $(\mathcal{A}^T\mathcal{A})^{-1}\mathcal{A}^T$.
This pseudo-inverse can then be used to compute the unknown control points
$G_{01}$, $G_{10}$, $G_{11}$, $G_{12}$, $G_{21}$ of the quadrilateral from the control points of the triangle parameterizations $T^w$.
However, while this approach has a strong theoretical foundation, we will see in our numerical examples that the method tends to result in self-intersecting quadrilateral parameterizations, even for simple geometries.
For this reason, we propose alternative averaging methods that are more application-oriented.

\paragraph{Parameterization averaging} 

In this approach, wee compute the biquadratic mapping $G \in \mathbb{Q}^2_2 \times \mathbb{Q}^2_2$ by averaging the triangle parameterizations $T^w \in \mathbb{P}^2_2 \times \mathbb{P}^2_2$, for $2\in\{1,2,3,4\}$.

Let us denote by $J_2[F](s_0,t_0)$ the $2$-jet, that is, the collection of $C^2$-data of the function $F$, evaluated at the point $(s_0,t_0)$. We then compute a biquadratic mapping that yields the best approximation to the $2$-jets of the four triangle parameterizations at the respective corner points, i.e.,
\begin{linenomath}\begin{equation*}
 \|J_2[G-G^1](0,0)\|^2+\|J_2[G-G^2](1,0)\|^2 +\|J_2[G-G^3](0,1)\|^2 + \|J_2[G-G^4](1,1)\|^2 \rightarrow \min.
\end{equation*}\end{linenomath}

This can be computed in terms of control points the following way: Let $G$ be represented in the biquadratic tensor-product Bernstein basis with control points $G_{ij}$
and let the mappings $T^w$ be represented in the same basis, with control points $T^w_{ij}$, for $i,j\in\{0,1,2\}$. The approximation of $2$-jets then results in
\[
 \left(\begin{array}{ccc}
  G_{02} & G_{12} & G_{22}  \\
  G_{01} & G_{11} & G_{21}  \\
  G_{00} & G_{10} & G_{20}
 \end{array}\right) =
 \left(\begin{array}{ccc}
  \frac13(T^{1}_{02}+T^{3}_{02}+T^{4}_{02}) & \frac12(T^{3}_{12}+T^{4}_{12}) & \frac13(T^{2}_{22}+T^{3}_{22}+T^{4}_{22}) \\
  \frac12(T^{1}_{01}+T^{3}_{01}) & \frac14(T^{1}_{11}+T^{2}_{11}+T^{3}_{11}+T^{4}_{11}) & \frac12(T^{2}_{21}+T^{4}_{21}) \\
  \frac13(T^{1}_{00}+T^{2}_{00}+T^{3}_{00}) & \frac12(T^{1}_{10}+T^{2}_{10}) & \frac13(T^{1}_{20}+T^{2}_{20}+T^{4}_{20})
 \end{array}\right).
\]

In terms of the coefficients $T^w_{ijk}$ of the triangle parameterizations, this results in
\begin{align*}
 G_{12} &= \frac{1}{2}(T^3_{011} + T^4_{011}),  \quad G_{01}  = \frac{1}{2}(T^1_{101} + T^3_{101}),\quad
 G_{21} = \frac{1}{2}(T^2_{011} + T^4_{110}),
  \quad G_{10} = \frac{1}{2}(T^1_{110} + T^2_{110}),\\
  G_{11} &= \frac{1}{8}\left(\sum_{w=1}^4 (T^w_{110} + T^w_{101} + T^w_{011}) - G_{00} - G_{02} - G_{20} - G_{22}\right).
\end{align*}

\paragraph{Control point averaging} 
As a very simple alternative to the previous approaches, we compute the control points of $G$ by directly averaging the control points of the triangles as follows
\[
 \left(\begin{array}{ccc}
  G_{02} & G_{12} & G_{22}  \\
  G_{01} & G_{11} & G_{21}  \\
  G_{00} & G_{10} & G_{20}
 \end{array}\right) =
 \left(\begin{array}{ccc}
  G_{02} & \frac{1}{2}(T^3_{011} + T^4_{011}) & G_{22} \\
  \frac{1}{2}(T^1_{101} + T^3_{101}) & \frac{1}{4}(T^1_{011} + T^2_{101} + T^3_{110} + T^4_{101}) & \frac{1}{2}(T^2_{011} + T^4_{110}) \\
  G_{00} & \frac{1}{2}(T^1_{110} + T^2_{110}) & G_{20}
 \end{array}\right).
\]
This is a more heuristic approach compared to Variants 1 and 2, but the obtained results are nonetheless useful and, in some cases, even better than the ones obtained using the other methods.

\subsubsection{Multi-patch domain}
For multi-patch domains, we assume that the computational domain is segmented into quadrilateral patches in a $C^0$-conforming way, without hanging vertices. We then perform the procedure presented in the previous section to compute biquadratic reparameterizations for all subdomains independently. To obtain a $C^0$-conforming discretization we then average the corresponding control points on all interfaces.

\subsection{Sampling strategies}\label{sec:point-sampling}
Depending on the modeling problem, we aim to minimize the best approximation error in the isogeometric space either with respect to the unknown solution \eqref{eq:L2-error} or the derivatives of the unknown solution \eqref{eq:L2-derivative-error}.

Our method minimizes the distance between the sampled points and the (unknown) graph surface:
\begin{linenomath}\begin{equation}\label{eq:graph-sampling-error}
 \sum_i \left|(x_i,y_i,z_i)^T-(G_x(s_i,t_i),G_y(s_i,t_i),\hat{u}_h(s_i,t_i))^T\right|^2 \rightarrow \min_{G,\hat{u}_h,s_i,t_i}.
\end{equation}\end{linenomath}

In order to optimize the approximation error with respect to the correct norm, we sample the input data for our method problem-dependent either from the graph surface of an initial approximation to the solution or from the graph surface of the derivatives of the initial approximation.

To this end, we first compute the initial solution
$u_{\mathit{init}} := \Pi^{\mathit{init}}_{\theta} (u)$ using the initial parameterization $G_{\mathit{init}}$ of the computational domain $\Omega$ and an initial mesh size $\theta$.

In order to achieve scaling invariance of the result of the standardization of the neural network input, see Section~\ref{NN}, with respect to the function
as well as the size of the computational domain
we replace $u_{\mathit{init}}$ by the normalized function
\begin{linenomath}\begin{equation*}
\frac{\max_k\left(\mathrm{diam}\,\Omega^k\right)}{\max_\Omega(u_{\mathit{init}}) - \min_\Omega(u_{\mathit{init}})}u_{\mathit{init}}.
\end{equation*}\end{linenomath}
\subsubsection{$L^2$-projection}
In the case of $L^2$-projection~\eqref{eq:L2-error}, we then sample values directly from the graph surface of $u_{\mathit{init}}$, i.e.,
$z_i:=u_{\mathit{init}}(x_i,y_i) = \hat u_{\mathit{init}}(s_i,t_i)$,
where $(s_i, t_i) = G^{-1}_{\mathit{init}}(x_i, y_i)$.
This yields a similar result to minimizing the $L^2$-error, which can be approximated by sampling points:
\[
 \| ({\rm I}-\Pi_{L^2,h}) (u) \|_{L^2(\Omega)} \sim \sum_i \left|u(x_i,y_i)-u_h(x_i,y_i)\right|^2.
\]

\subsubsection{$H^1$-projection}\label{sec:samplingh1}
In the case of $H^1$-projection~\eqref{eq:L2-derivative-error}, we additionally sample values from different directional derivatives of $\hat u_{\mathit{init}}$, i.e., of the spline function defined on the parameter domain $\hat\Omega$. More precisely, we sample values
$z_i:= \frac{\partial}{\partial \vec v} \hat u_{\mathit{init}}(s_i,t_i) = \frac{\partial}{\partial \vec v}\left(u_{\mathit{init}} \circ G_{\mathit{init}}\right)(s_i,t_i)$,
where $\vec v\in\mathbb R^2$.
Since $u_{\mathit{init}}$ is in an isogeometric function space $V_{\theta, G_{\mathit{init}}}$, $u_{\mathit{init}} \circ G_{\mathit{init}}$ is a tensor-product B-spline function and therefore the directional derivatives are simple to evaluate.

In particular, if we choose $\vec v$ to be $(\frac1{\sqrt2},\frac1{\sqrt2})$ or $(\frac1{\sqrt2},-\frac1{\sqrt2})$, we obtain the derivatives in direction of the corners of the quadrilateral patches. If we choose $\vec v$ to be either $(1,0)$ or $(0,1)$, we obtain the derivatives in the direction of the edges of the quadrilateral patch. In our numerical experiments, we analyze the influence of choosing different directions $\vec v$.

Applying the method to sets of points sampled from the initial solution as well as different directional derivatives leads to several different reparameterizations $G^r$ of the computational domain. In order to obtain the final reparameterization, we choose for each control point $G_{ij}$ the control point $G_{ij}^r$ with the largest distance to the corresponding control point $G_{\mathit{init},ij}$ of the initial parameterization. This ensures that we capture the deformation of the parameterization in all regions where an increased mesh density is needed.

\section{Numerical experiments}\label{sec:numerical-experiments}
In this section, we present the results of our numerical experiments.
On the one hand, we aim to demonstrate the suitability of our approach for optimizing parameterizations for solving PDEs using isogeometric analysis.
On the other hand, we want to motivate the choice of the averaging strategy, see Section~\ref{biquadratic}, and sampling strategy, see Section~\ref{sec:point-sampling}.
We start by comparing the different averaging approaches on a number of single patch domains. Afterwards we focus on multi-patch domains, comparing the different sampling strategies when solving the Poisson equation. Finally, we apply our method to some examples inspired by real-world problems.

The isogeometric methods for solving the PDE and $L^2$-projection was implemented in C\texttt{++} using the open source library G+Smo, cf.~\cite{gismoweb}. Our method for optimizing the isogeometric function spaces was implemented using Python, with the neural network being defined using PyTorch, cf.~\cite{Paszke}.
As batch size when optimizing the triangular B\'ezier parameterizations we set $M=20$ and we stopped the optimization of each triangle after the $L^2$-error stops improving for two steps in a row.

\subsection{Comparison of the different averaging approaches}
We apply our method to different computational domains both for the $L^2$-projection and the Poisson problem and compare the three averaging approaches described in Section~\ref{biquadratic}: Best-approximating inverse (BAI), control point averaging (CPA) and parameterization averaging (PA).

As an illustration of the overall process, we start by considering the $L^2$-projection problem with exact solution $u(x,y) = ((x-1)^2+(y-1)^2)^{\frac{1}{16}}$ over the physical domain $\Omega = [0,1]^2$.
First of all, $\Omega$ is divided into four triangles as described in Figure~\ref{fig:Quad_Splitting}.

\begin{figure}
	\begin{center}
		\subfloat[$\Delta_1$]{\includegraphics[scale=0.3, viewport=3.5cm 1cm 13cm 11cm, clip]{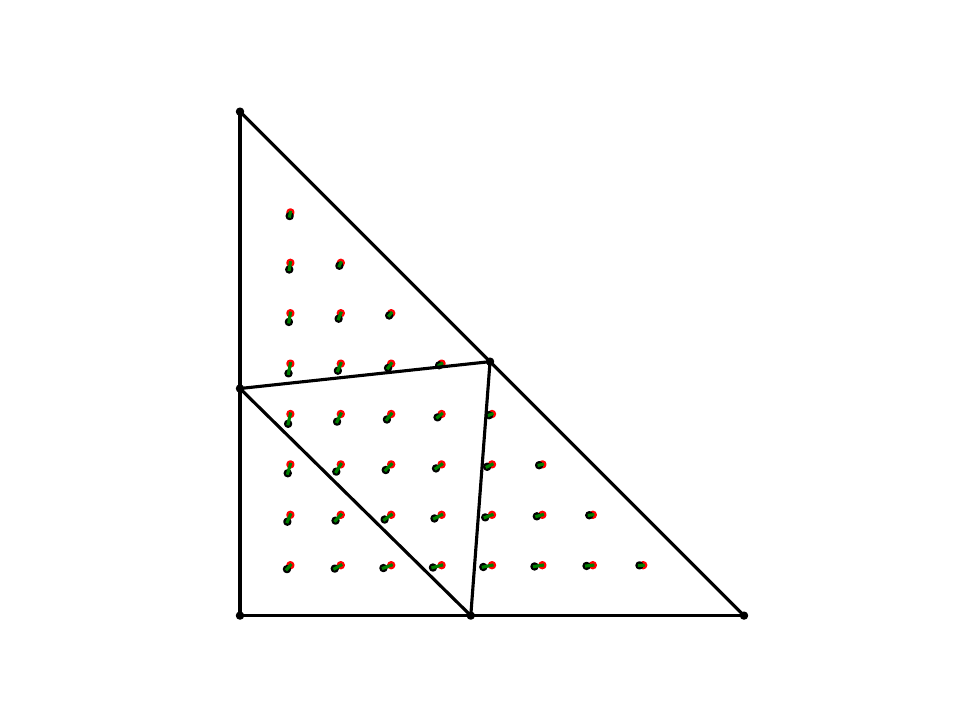}}
		\subfloat[$\Delta_2$]{\includegraphics[scale=0.3, viewport=3.5cm 1cm 13cm 11cm, clip]{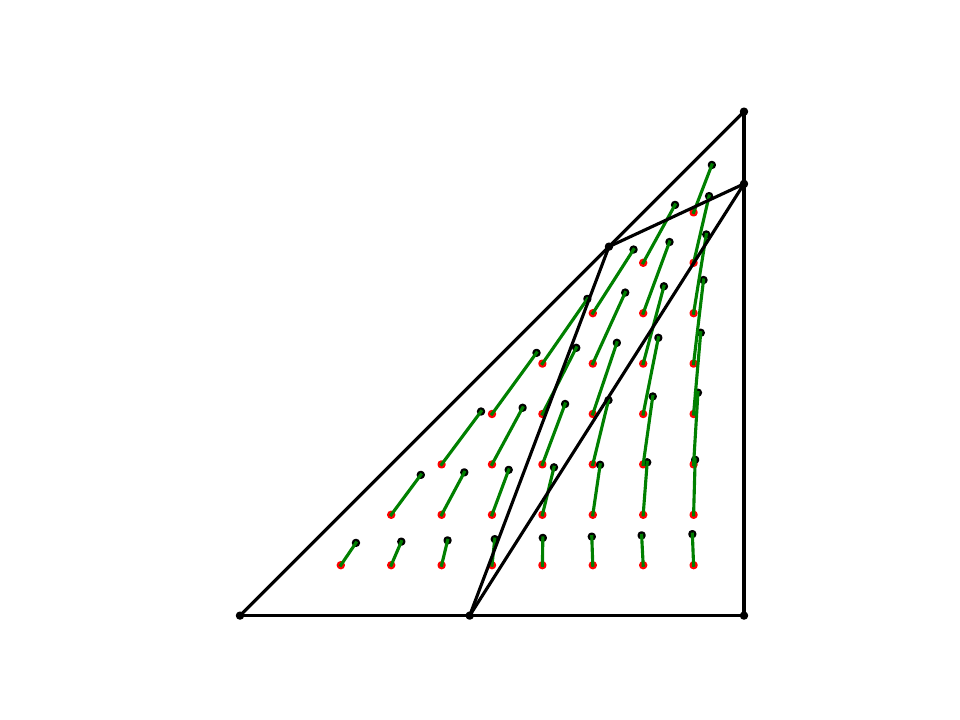}}
		\subfloat[$\Delta_3$]{\includegraphics[scale=0.3, viewport=3.5cm 1cm 13cm 11cm, clip]{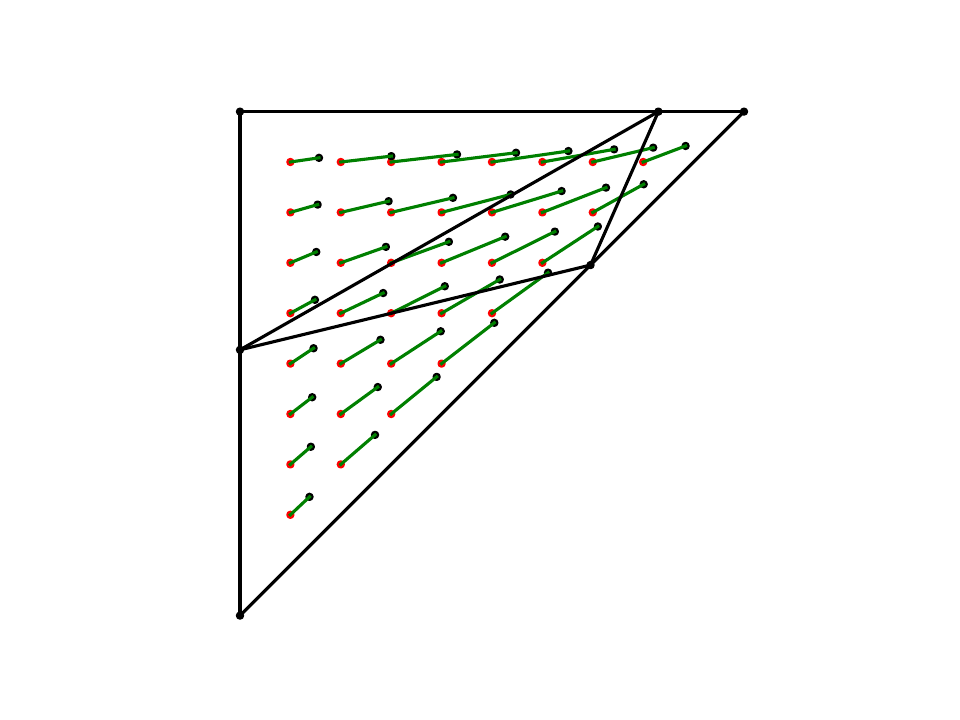}}
		\subfloat[$\Delta_4$]{\includegraphics[scale=0.3, viewport=3.5cm 1cm 13cm 11cm, clip]{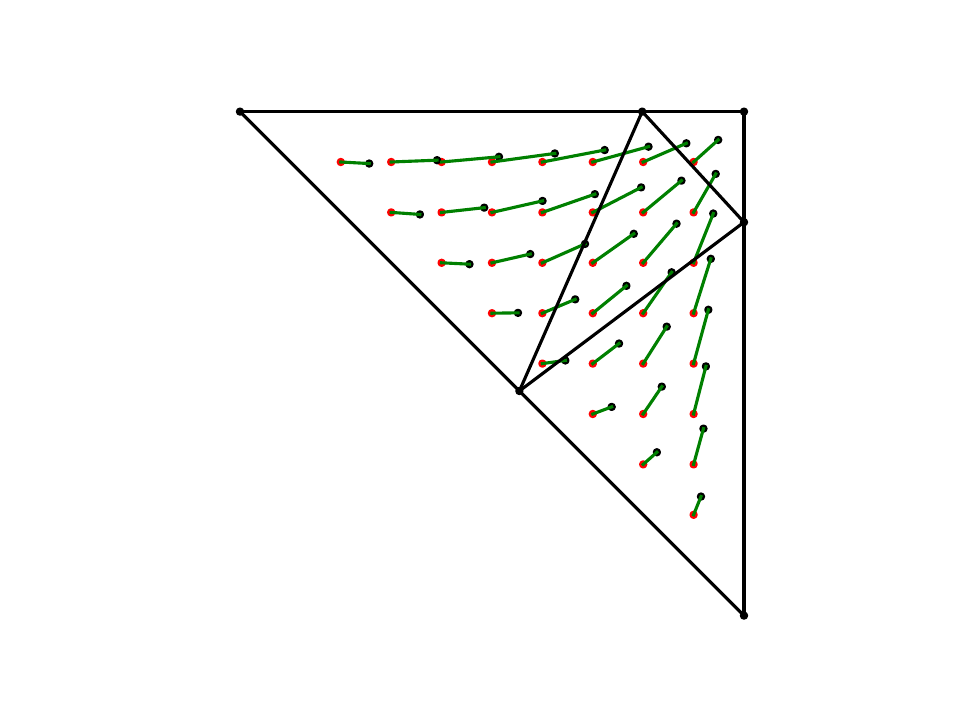}}\\
	\end{center}
	\vspace{-0.5cm}
	\caption{Triangle deformations. The deformed triangular B\'ezier control nets are shown in black. The green lines visualize the point-wise change of the parameterization from undeformed (red) to deformed (black).}\label{triangletoquad}
\end{figure}

\begin{figure}
	\centering
	\subfloat[Best-approximating inverse]{\includegraphics[scale=0.3, viewport=0cm 1cm 15.5cm 10.5cm, clip]{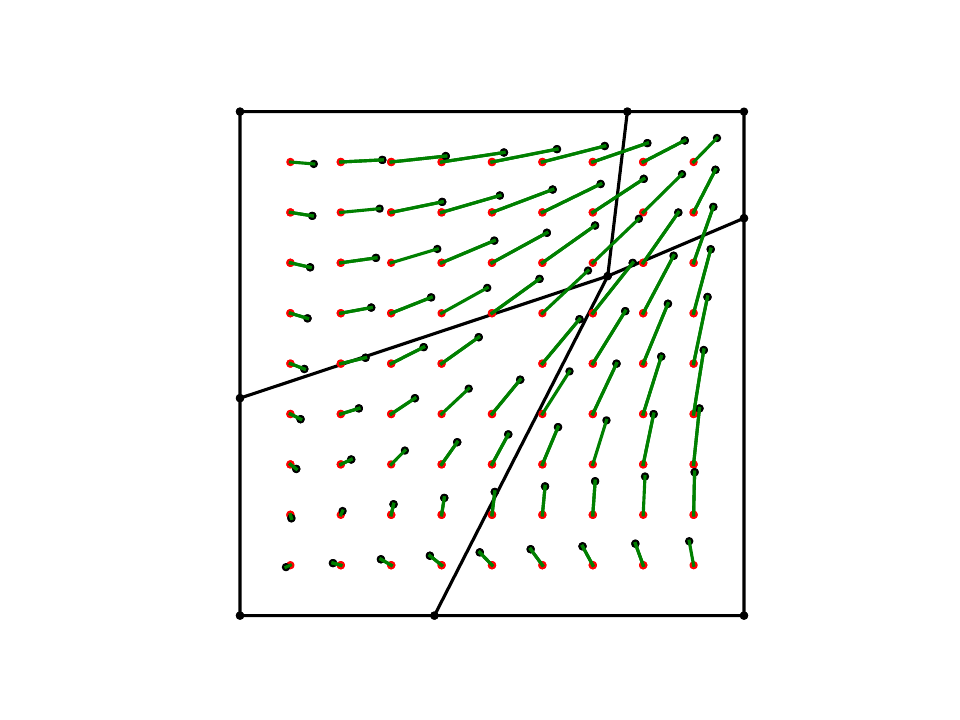}}
	\subfloat[Parameterization averaging]{\includegraphics[scale=0.3, viewport=0cm 1cm 15.5cm 10.5cm, clip]{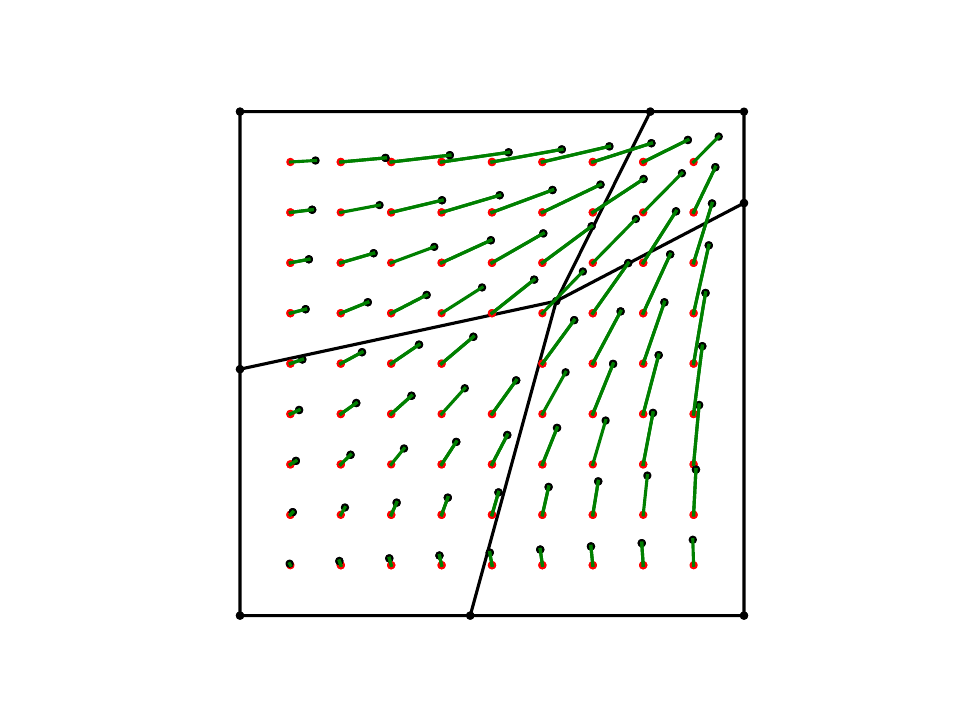}}
	\subfloat[Control point averaging]{\includegraphics[scale=0.3, viewport=0cm 1cm 15.5cm 10.5cm, clip]{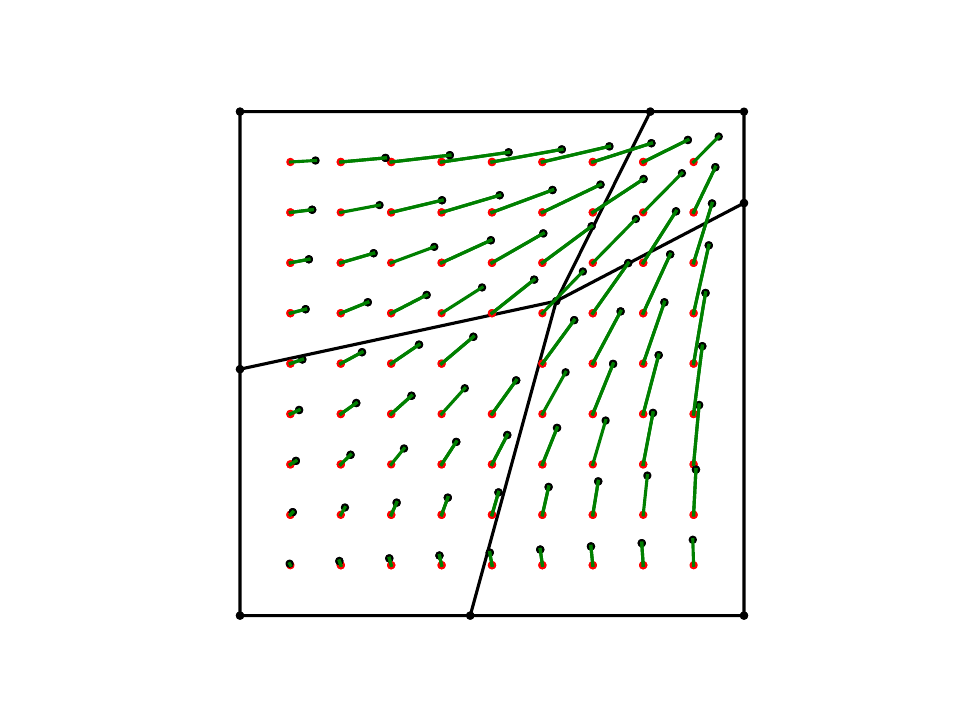}}
	\caption{Comparison of the averaging approaches presented in Section~\ref{biquadratic}. The deformed tensor-product B\'ezier control nets are shown in black. The green lines visualize the point-wise change of the parameterization from undeformed (red) to deformed (black).}
	\label{comparison_averaging}
\end{figure}

Then we apply the method described in Section~\ref{sec:triangle} to each triangular domain individually. Figure~\ref{triangletoquad} shows the resulting deformed B\'ezier nets for each triangle.
One can observe that the control points of the reparameterization of $\Delta_2, \Delta_3$ and $\Delta_4$ move closer to the singularity (top-right corner). 
On the other hand, one observe only a slight deformation in the reparameterization of $\Delta_1$ since the singularity was not present in that part of the domain.
Figure~\ref{comparison_averaging} shows the three quadrilateral control nets resulting from applying the three averaging approaches to these triangular B\'ezier nets. 

In the next examples, we consider different domains and we analyze the behavior of the $L^2$-error and the $H^1$-error when solving the $L^2$-projection problem and the Poisson problem.

\begin{figure}
	\centering
	\subfloat[Exact solution]{\includegraphics[scale=0.18, viewport=8cm 1cm 32.9cm 23cm, clip]{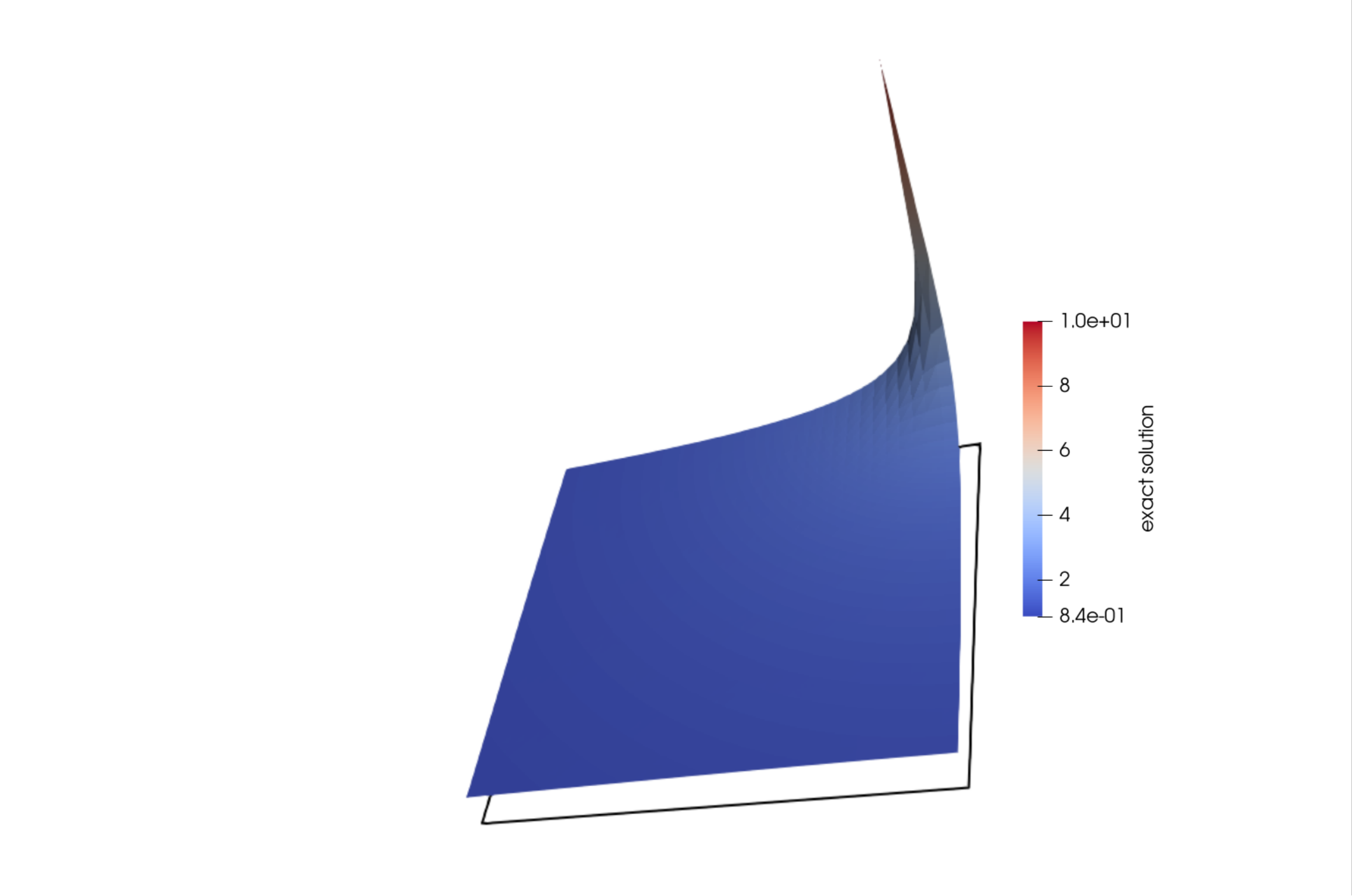}\label{L2_SingeSquareSol}} \hspace{2cm}
	\subfloat[$L^2$-errors]{
		\includegraphics{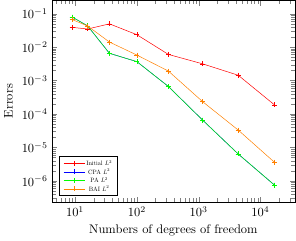}
		\label{L2_SingeSquareErrors}
	}
	\caption{$L^2$-projection on a square domain with a corner singularity: Exact solution and convergence of $L^2$-errors.}
\end{figure}

\paragraph{Square domain with a corner singularity}
We consider the $L^2$-projection problem with exact solution $u(x,y) = ((x - 1)^ 2 + (y - 1)^2 + 0.0001)^{-\frac{1}{4}}$ defined over the unit square $[0,1]^2$, Figure~\ref{L2_SingeSquareSol}.
While this function is in $C^\infty$ (and in $H^k$ for any $k\in\mathbb{N}$), the gradient becomes very steep close to $(1,1)$. The unperturbed function is in $L^2$, but not in $H^1$.

\begin{figure}
	\centering 
	\subfloat[Best-approximating inverse]{\includegraphics[scale=0.2, viewport=0cm 0cm 13.5cm 14.5cm, clip]{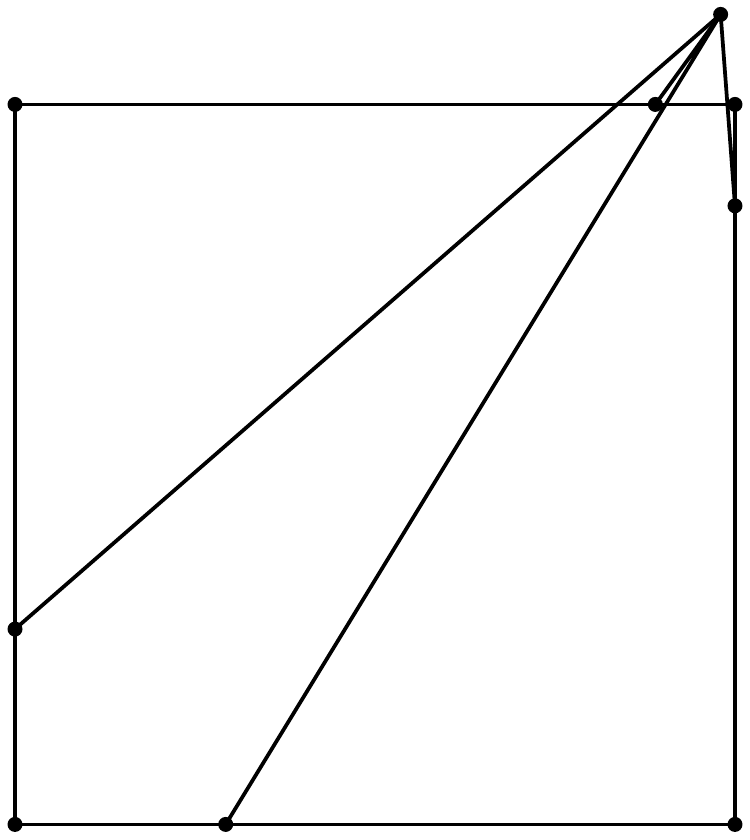}\includegraphics[scale=0.21, viewport=0cm 0cm 13cm 13cm, clip]{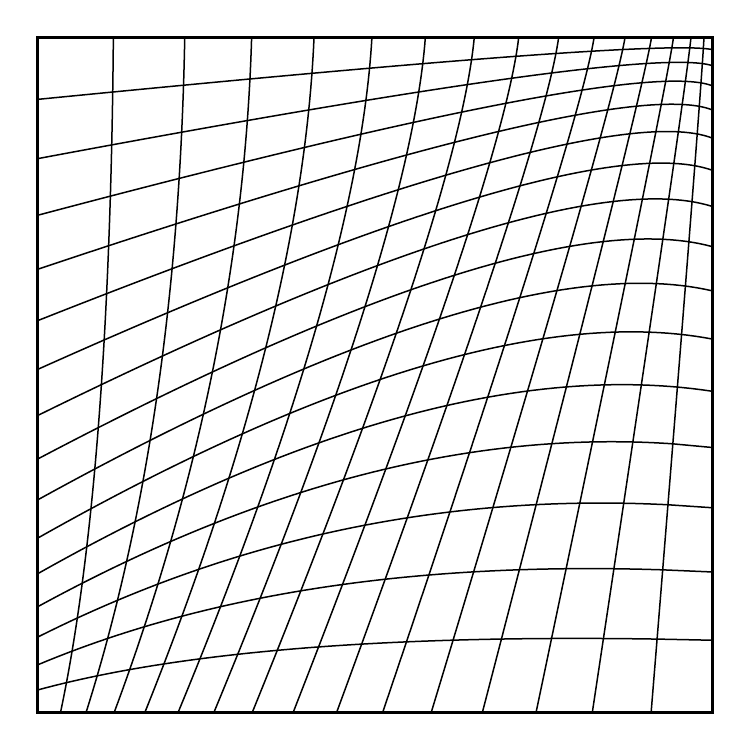}}
	\subfloat[Parameterization averaging]{\includegraphics[scale=0.2, viewport=0cm 0cm 13cm 13cm, clip]{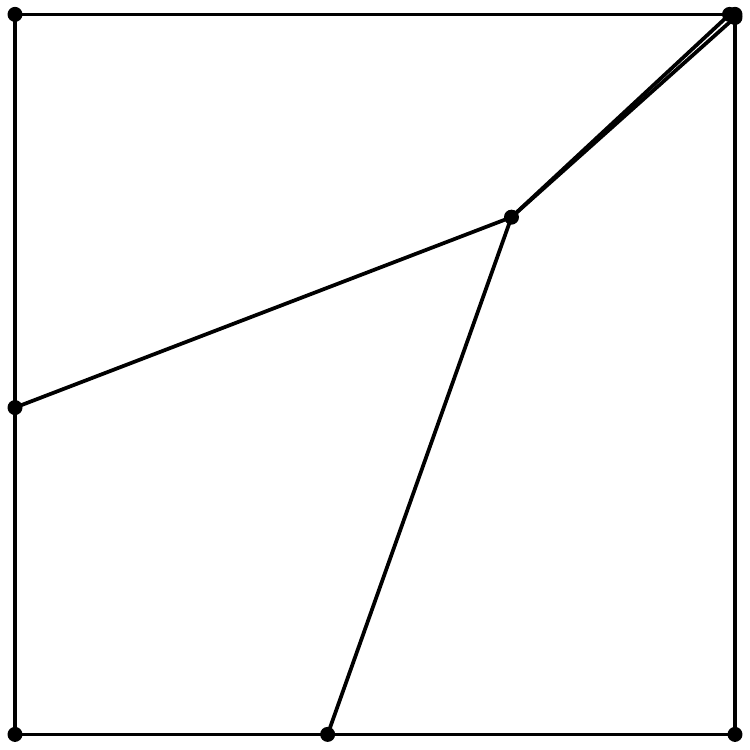}\includegraphics[scale=0.21, viewport=0cm 0cm 13cm 13cm, clip]{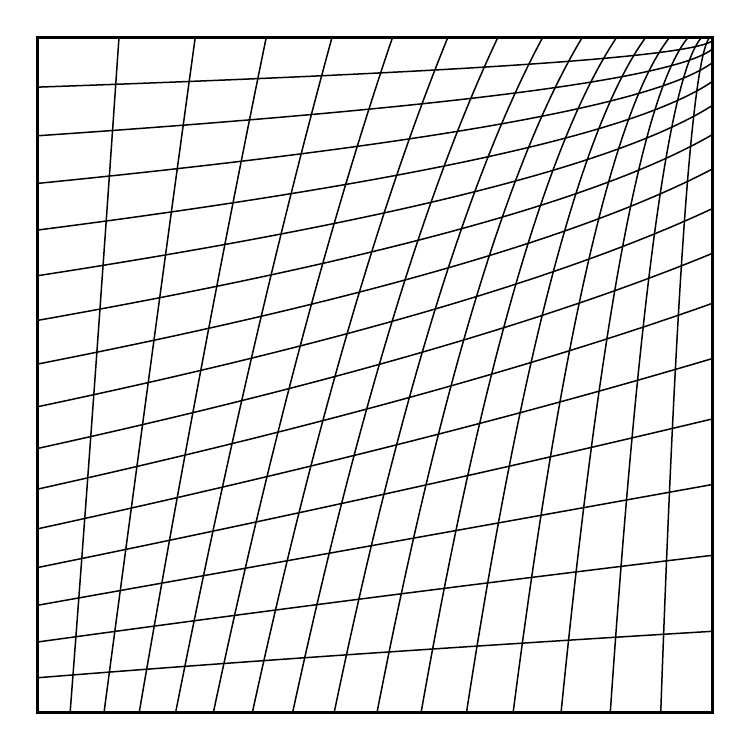}}
	\subfloat[Control point averaging]{\includegraphics[scale=0.2, viewport=0cm 0cm 13.5cm 13cm, clip]{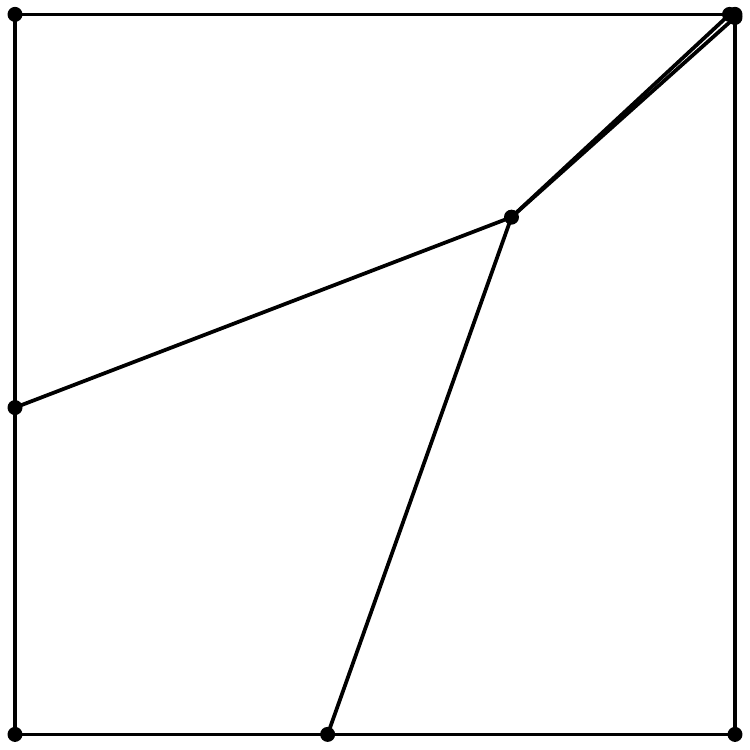}\includegraphics[scale=0.21, viewport=0cm 0cm 13cm 13cm, clip]{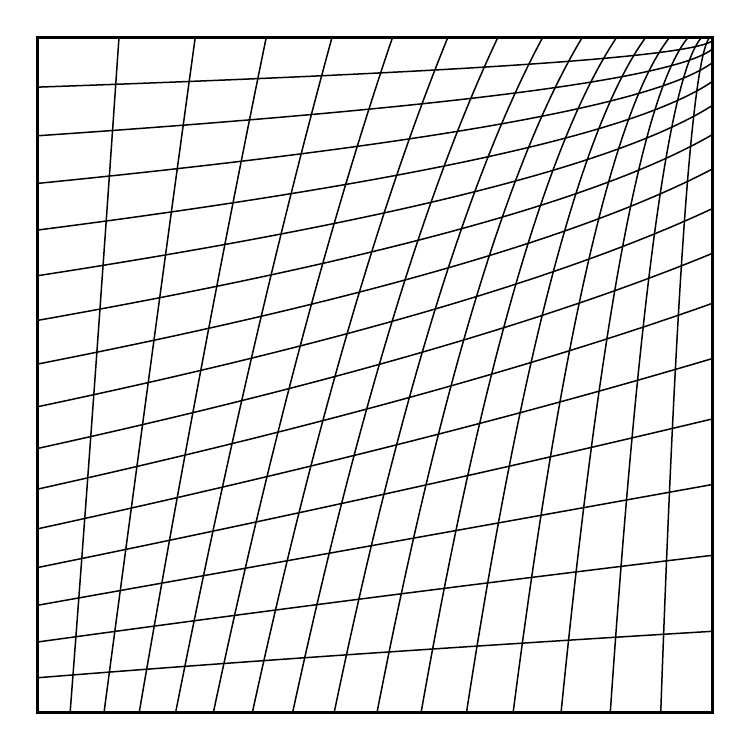}}
	\caption{$L^2$-projection on a square domain with a corner singularity:  Deformed control nets for the three averaging approaches.}
	\label{L2_SingeSquareDeformCP}
\end{figure}

In Figure~\ref{L2_SingeSquareDeformCP} we present the deformed control nets and the parameter lines of the physical domain for each averaging approach. We observe that for all three approaches, the control points move closer to $(1,1)$. 
 Figure~\ref{L2_SingeSquareErrors} presents the $L^2$-errors corresponding to the original and deformed control nets, for each averaging approach, after applying several iterations of $h$-refinement.
  One can observe an improvement of the $L^2$-errors of more than two orders of magnitude  when using the deformed control nets of the 
  control point averaging and the parameterization averaging approaches, while the best-approximating inverse approach improves the error by more than one order of magnitude.


\begin{figure}
	\centering
	\subfloat[Exact solution]{\includegraphics[scale=0.18, viewport=8cm 0cm 31.2cm 22cm, clip]{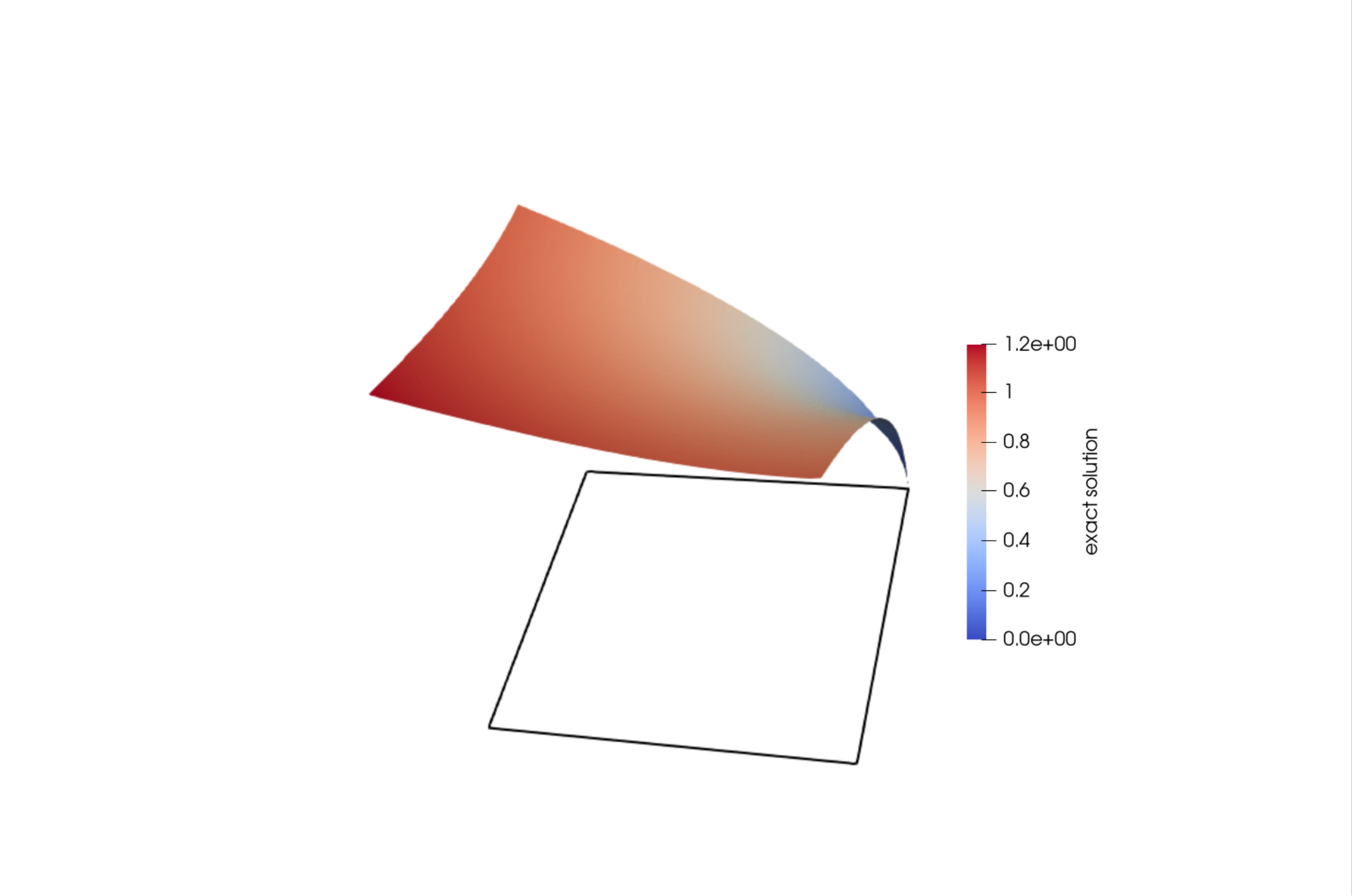}\label{Po_SingleSquareSol}} \hspace{2cm}
	\subfloat[$L^2$- and $H^1$-errors]{
		\includegraphics{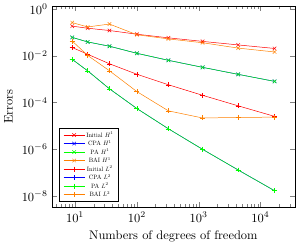}\label{Po_SingleSquareErrors}
	}
	\caption{Poisson problem on a square domain with a corner singularity: Exact solution and convergence of the errors.}
\end{figure}

On the same domain, we now consider the Poisson problem with Dirichlet boundary conditions, as described in Section~\ref{sec:modelProblems}.
We choose the right-hand side $f$ as well as the boundary function $g_D$ according to the exact solution $u(x,y)=((x - 1)^ 2 + (y - 1)^ 2)^\frac{1}{4}$ , see Figure~\ref{Po_SingleSquareSol}. 

\begin{figure}
	\centering
	\subfloat[Best-approximating inverse]{\includegraphics[scale=0.215, viewport=0cm 0cm 13cm 15cm, clip]{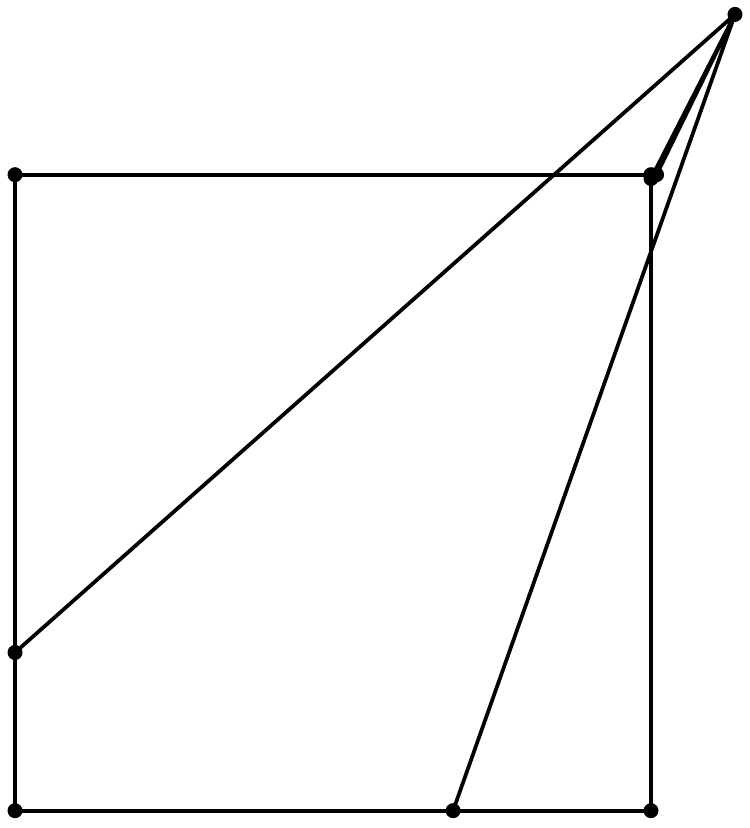}\includegraphics[scale=0.208, viewport=0cm 0cm 13cm 15cm, clip]{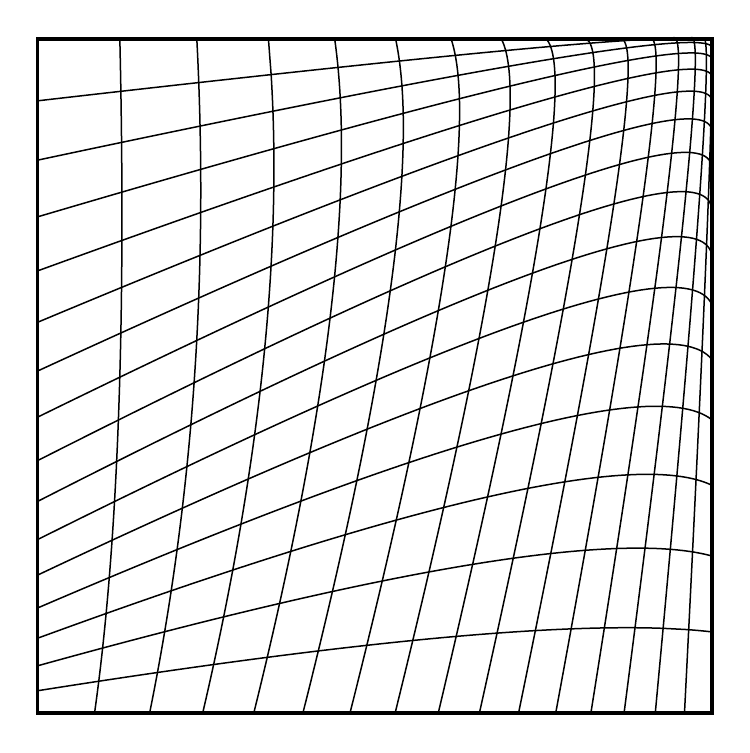}}
	\subfloat[Parameterization averaging]{\includegraphics[scale=0.2, viewport=0cm 0cm 13cm 15cm, clip]{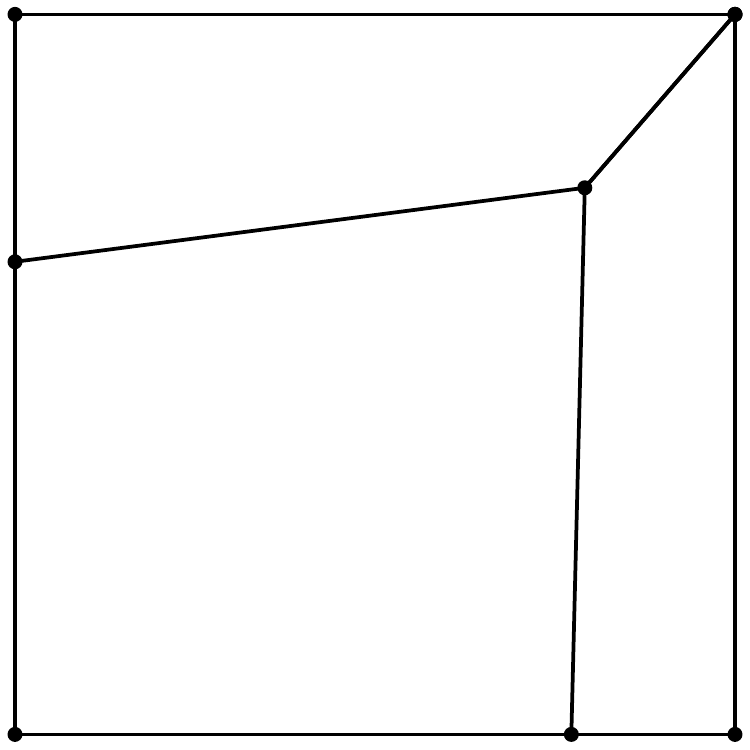}\includegraphics[scale=0.208, viewport=0cm 0cm 13cm 15cm, clip]{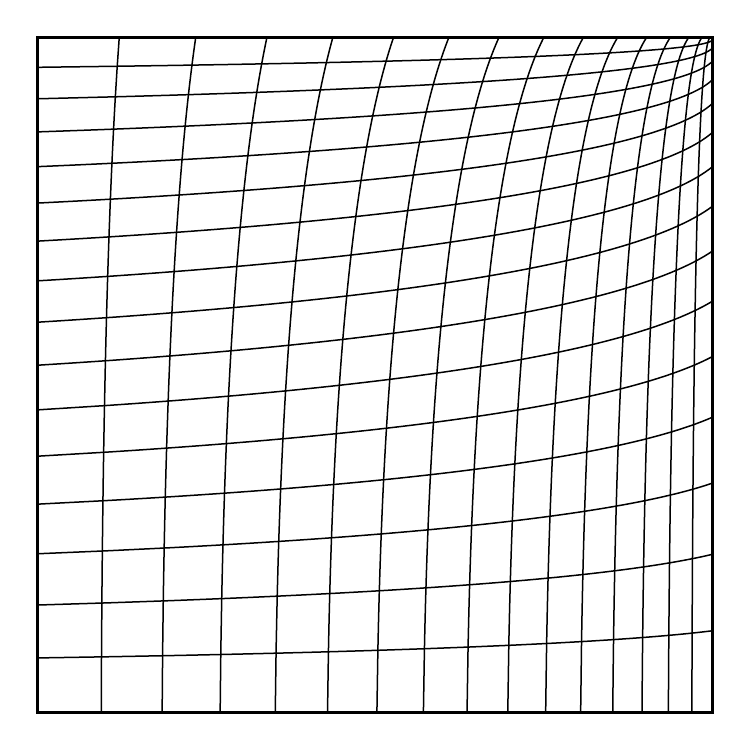}}
	\subfloat[Control point averaging]{\includegraphics[scale=0.2, viewport=0cm 0cm 13.5cm 15cm, clip]{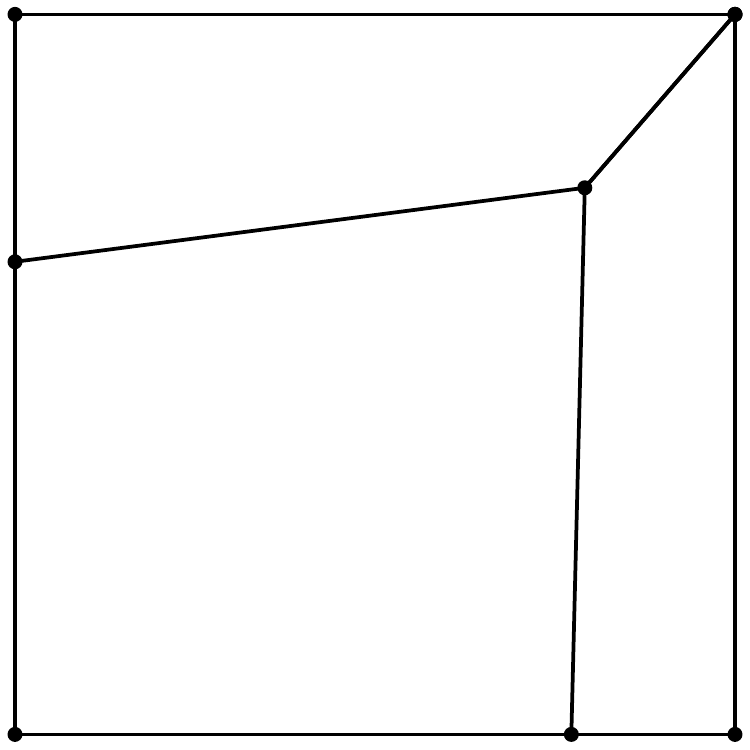}\includegraphics[scale=0.208, viewport=0cm 0cm 13cm 15cm, clip]{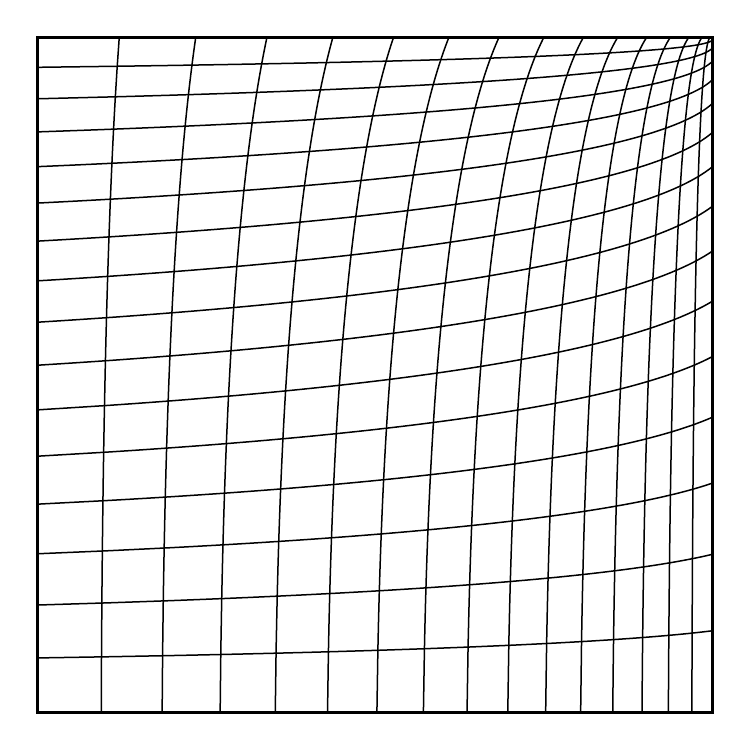}}
	\caption{Poisson problem on a square domain with a corner singularity:  Deformed control nets for the three averaging approaches.}\label{Po_SingleSquareDeformCP}
\end{figure}

Figure~\ref{Po_SingleSquareDeformCP} shows the control meshes and parameter lines of the reparameterizations resulting from the three averaging approaches after applying our method.
While the control point averaging and the parameterization
averaging result in regular parameterizations of the domain, the best-approximating inverse
approach leads to a non-regular map since the central control point moves too far outside of the domain.

Figure~\ref{Po_SingleSquareErrors} show the error convergence for the three approaches.
For the control point averaging and the parameterization averaging we observe an improvement of more than one order of magnitude for the $H^1$-error and of more than three orders of magnitude for the $L^2$-error. Moreover, the error as well as the control nets for these two approaches are very similar.
The best-approximating inverse initially improves the $L^2$-error, but after some steps of refinement, the error deteriorates due to the non-regularity of the parameterization.


\paragraph{Square domain with  side singularity}

\begin{figure}
	\centering
	\subfloat[Exact solution]{\includegraphics[scale=0.175, viewport=10cm 0cm 31.4cm 23cm, clip]{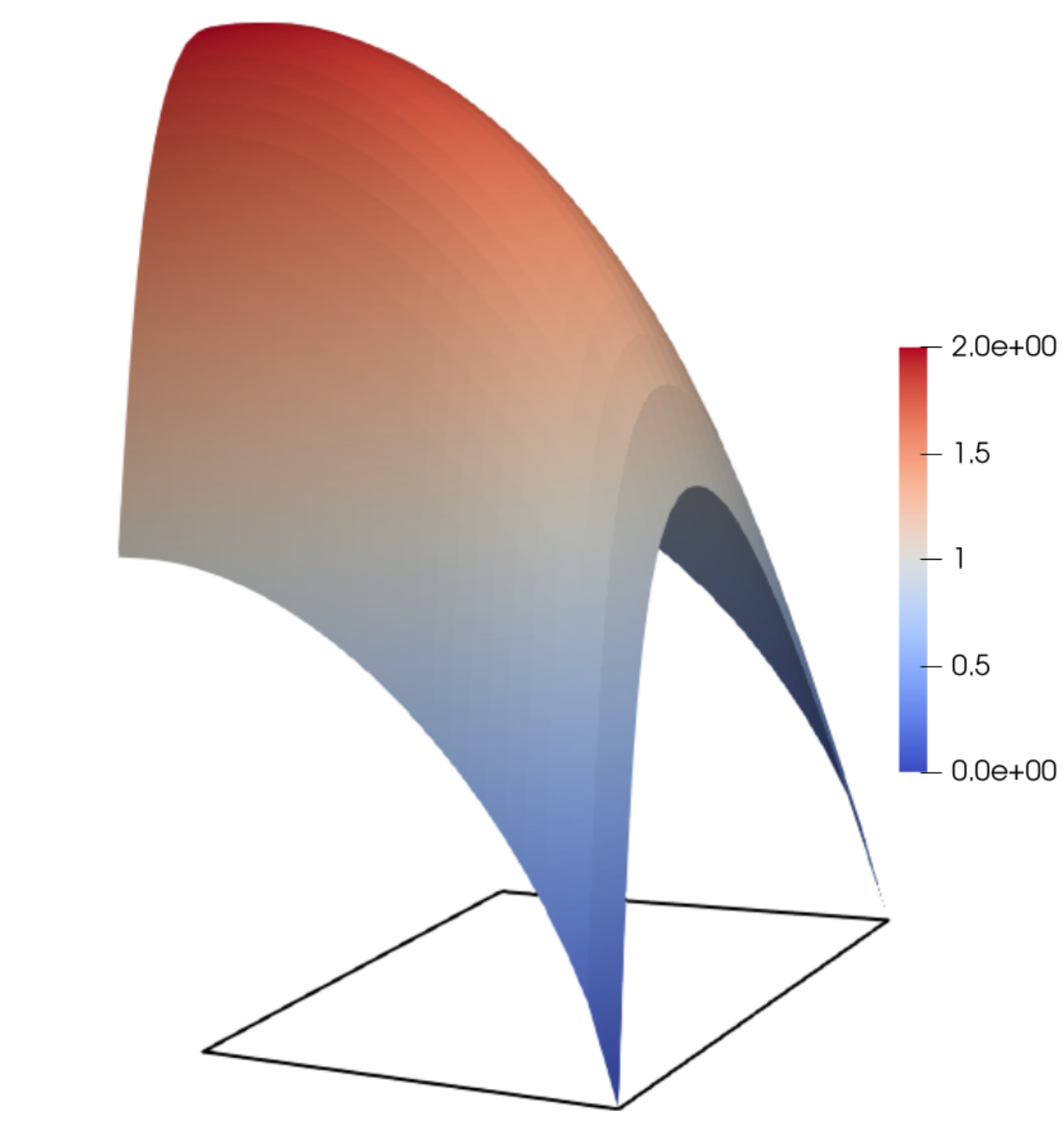}\label{Po_SingleSqureSideSingSol}} \hspace{2cm}
	\subfloat[$L^2$- and $H^1$-errors]{
		\includegraphics{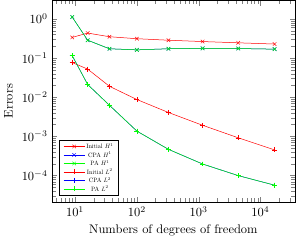} \label{Po_SingleSqureSideSingErrors}
	}
	\caption{Poisson problem on a square domain with side singularity: Exact solution and convergence.}
\end{figure}

In this example we solve the Poisson problem with Dirichlet boundary conditions and exact solution
$u(x,y)=\sin(\pi y) + (1 - x ^2)^{\frac35}$.
This function has a side singularity on the right boundary of the  physical domain $[0,1]^2$, see Figure~\ref{Po_SingleSqureSideSingSol}. 

\begin{figure}
	\centering
	\subfloat[Best-approximating inverse]{\includegraphics[scale=0.4, viewport=0cm 0cm 13cm 6.5cm, clip]{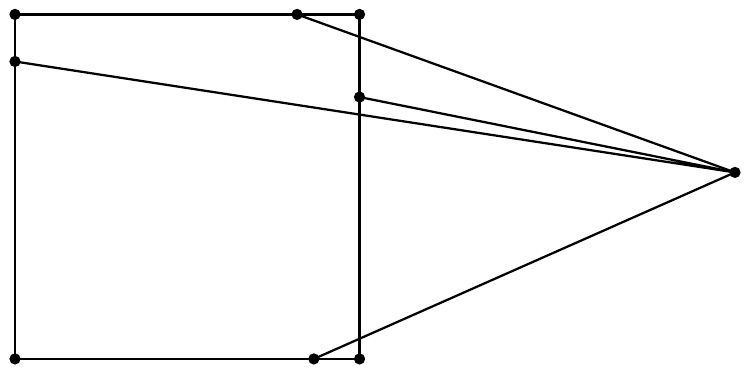}\includegraphics[scale=0.25, viewport=0cm 0cm 13cm 12cm, clip]{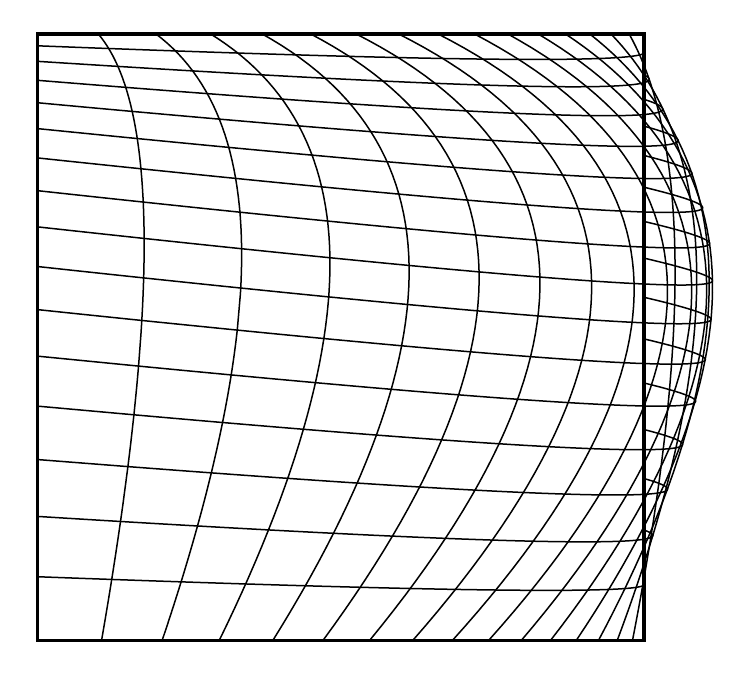}} \\
	\subfloat[Parameterization averaging]{\includegraphics[scale=0.2, viewport=0cm 0cm 13cm 13cm, clip]{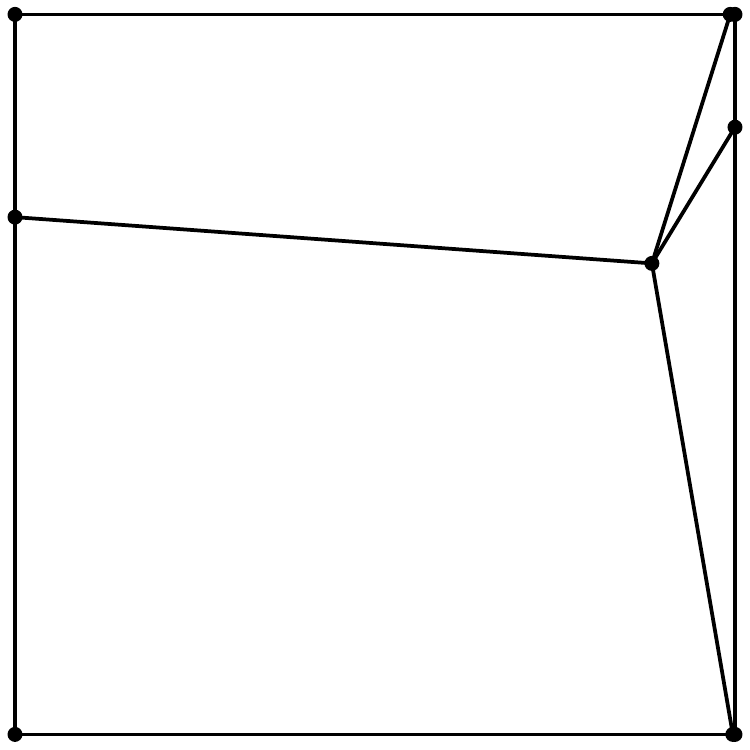}\includegraphics[scale=0.21, viewport=0cm 0cm 13cm 12.5cm, clip]{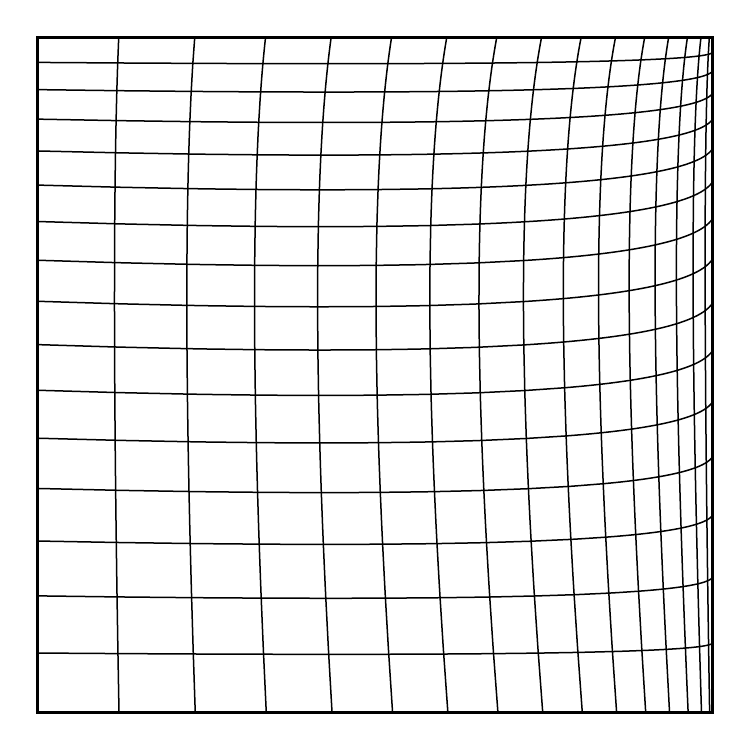}}
	\subfloat[Control point averaging]{\includegraphics[scale=0.2, viewport=0cm 0cm 13.5cm 13cm, clip]{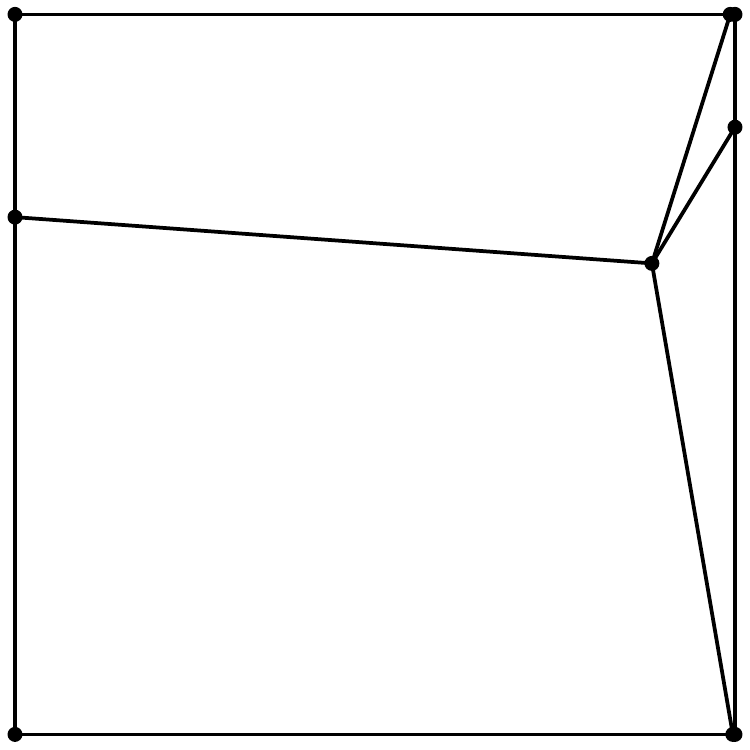}\includegraphics[scale=0.21, viewport=0cm 0cm 13cm 12.5cm, clip]{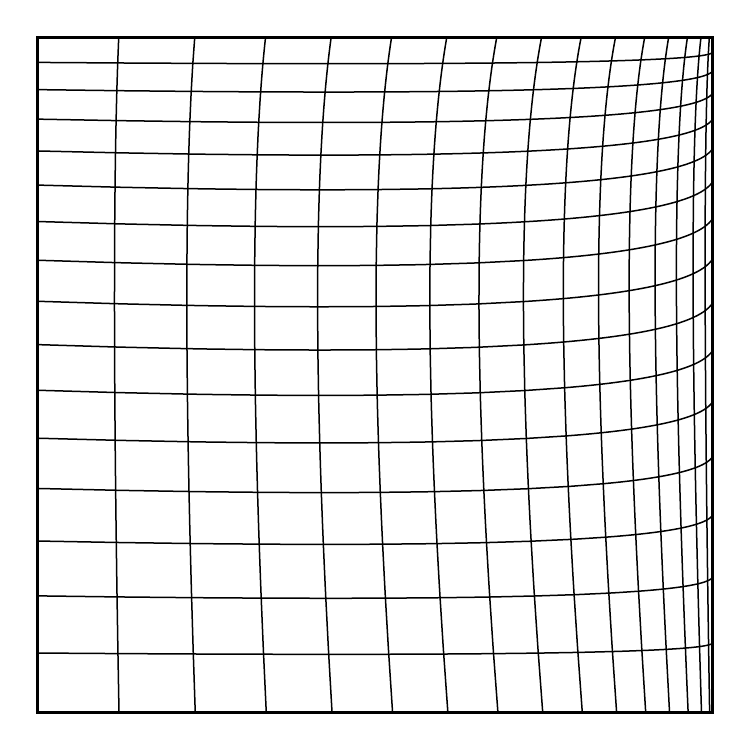}}
	\caption{Poisson problem on square domain with side singularity:  Deformed control nets for the three averaging approaches.}\label{Po_SingleSqureSideSingDeformCP}
\end{figure}

The control nets and the corresponding parameter lines are presented in Figure~\ref{Po_SingleSqureSideSingDeformCP} and
 Figure~\ref{Po_SingleSqureSideSingErrors} shows the $L^2$- and $H^1$-errors for the original and the deformed control nets after several steps of $h$-refinement.
 As in the previous example, the best-approximating inverse approach leads to a non-regular map that maps to points outside the computational domain. Since this is not a valid parameterization of the domain, we cannot solve the Poisson problem using this map, thus no errors are reported in Figure~\ref{Po_SingleSqureSideSingErrors}.
   For the remaining two averaging approaches, we observe an improvement of one order of magnitude in the $L^2$-error.
   The  $H^1$-error after two steps of refinement of the optimized parameterization is smaller than the $H^1$-error  after seven refinement steps of the original parameterization.
    We conclude that the best-approximating inverse approach is not suitable for our applications and will therefore not be considered in the following examples. Since the results for the parameterization averaging and the control point averaging are very similar in all examples, we from now on only perform parameterization averaging.


\subsection{Non-square quadrilateral domain with a corner singularity}
In the coming examples we consider a non-square quadrilateral domain $\Omega$.
 We consider the $L^2$-projection problem as well the Poisson problem. 
 
\begin{figure}
	\centering
	\subfloat[Exact solution]{\includegraphics[scale=0.22, viewport=14cm 5cm 30.5cm 25cm, clip]{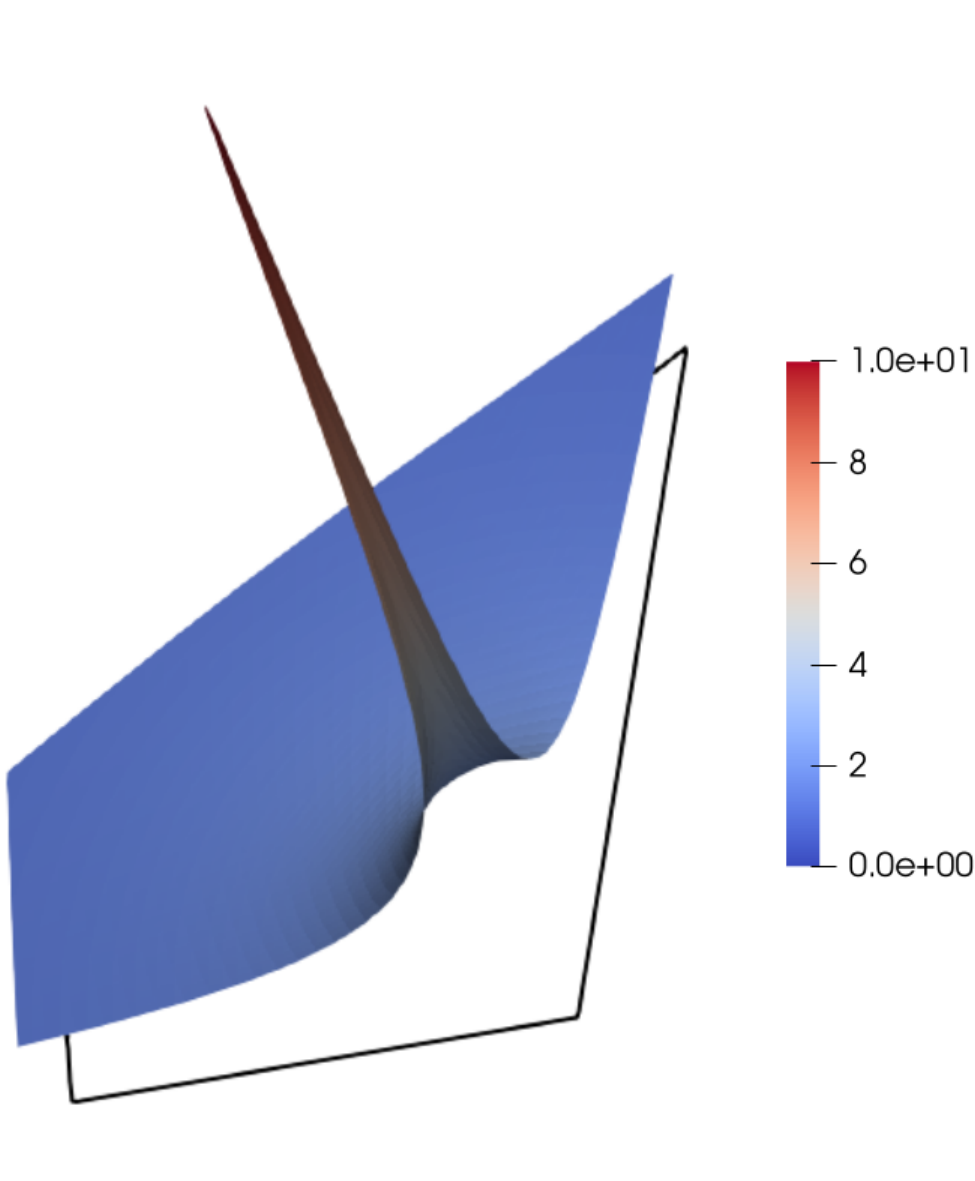}\label{L2_SingleQuadSol}}
	\subfloat[Resulting control net and  parameter lines]{\includegraphics[scale=0.28, viewport=-1cm -2cm 13cm 13cm, clip]{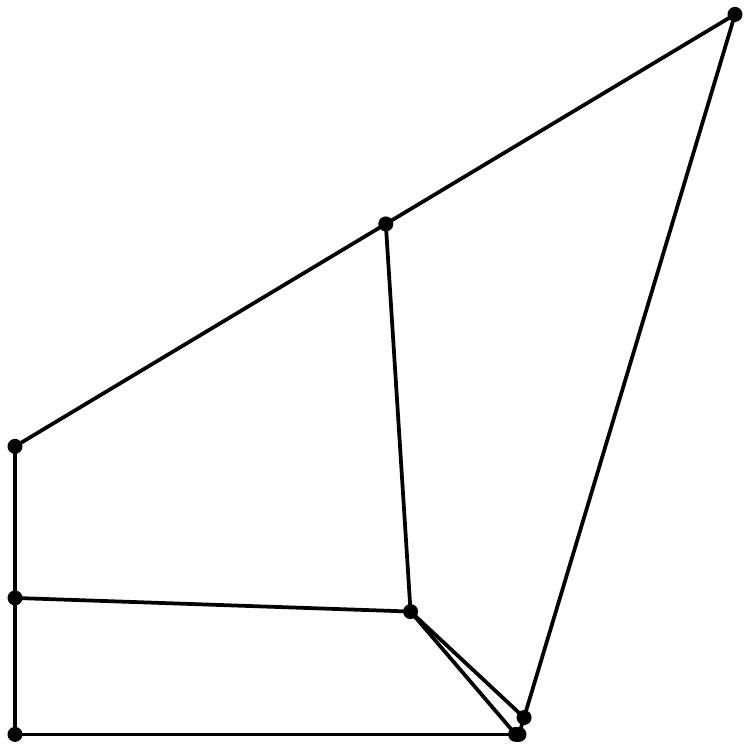}\includegraphics[scale=0.28, viewport=0cm -2cm 13cm 13cm, clip]{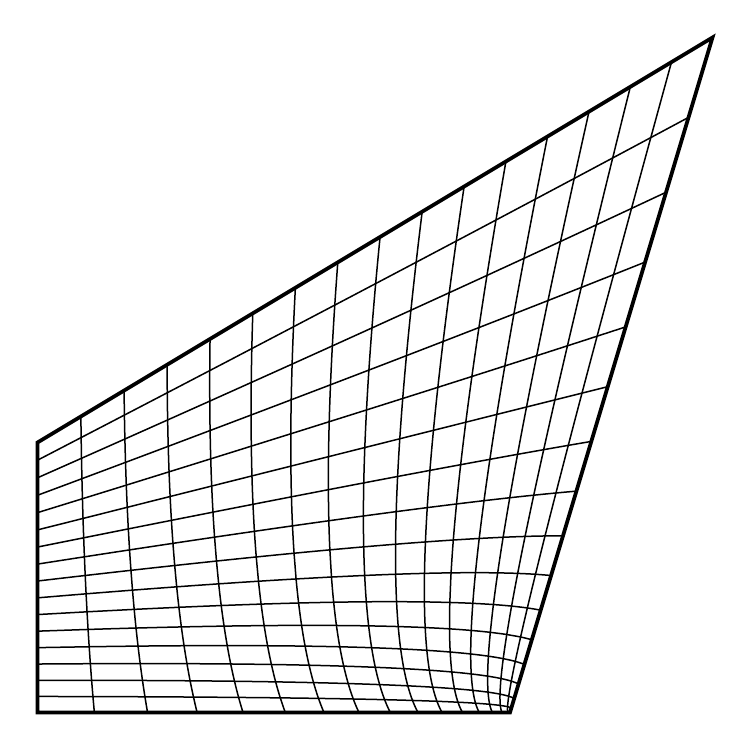}\label{L2_SingleQuadDeformCP}}
	\subfloat[$L^2$-error]{
		\includegraphics{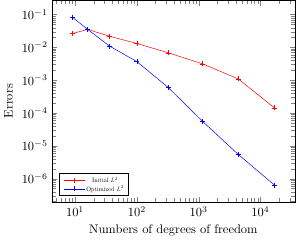}\label{L2_SingleQuadErrors}
	} 
	
	\caption{$L^2$-projection on a non-square quadrilateral domain with a corner singularity.}
\end{figure}

For the $L^2$-projection we consider the exact solution $u(x,y) = ((x - 0.7)^2 + y^2 + 0.0001)^{-\frac{1}{4}}$, see Figure~\ref{L2_SingleQuadSol}. Similar to the example in Figure~\ref{L2_SingeSquareSol} we perturb an otherwise singular function to obtain a function that is $C^\infty$, but with very high $H^1$-norm. The point $(0.7,0)$, where the gradient attains its maximum norm is at the lower right corner of the domain.
The deformed control net resulting from the parameterization averaging is shown in Figure~\ref{L2_SingleQuadDeformCP}.
 We present the $L^2$-error after several iterations of $h$-refinement in Figure~\ref{L2_SingleQuadErrors}. An improvement of more than two orders of magnitude in the $L^2$-error is observed.

\begin{figure}
	\centering
	\subfloat[Exact solution]{\includegraphics[scale=0.18, viewport=10cm -2cm 30.5cm 22cm, clip]{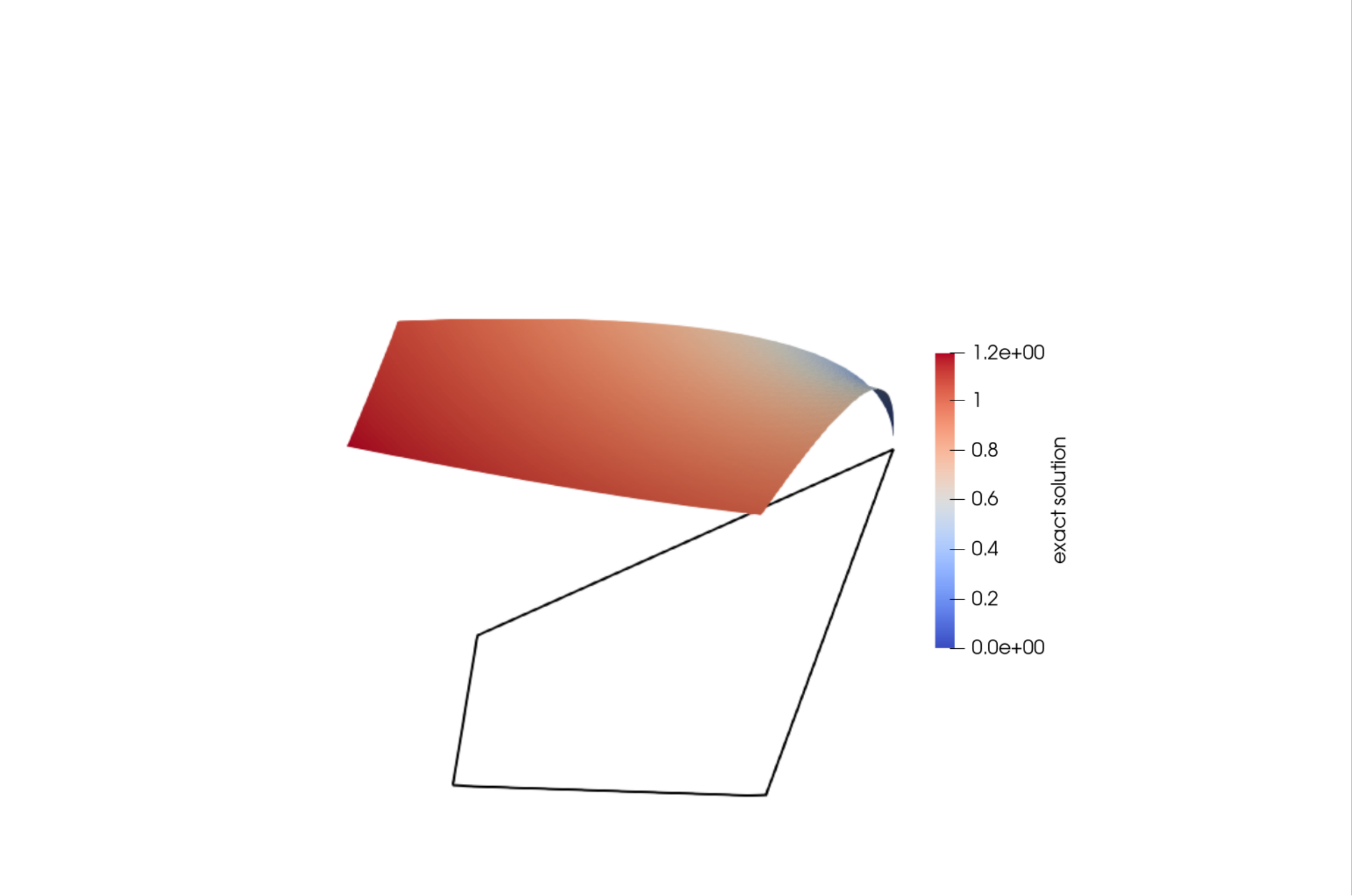}\label{Po_SingleQuadSol}}
	\subfloat[Resulting control net and  parameter lines]{\includegraphics[ scale=0.28, viewport=-1cm -2cm 13cm 13cm, clip]{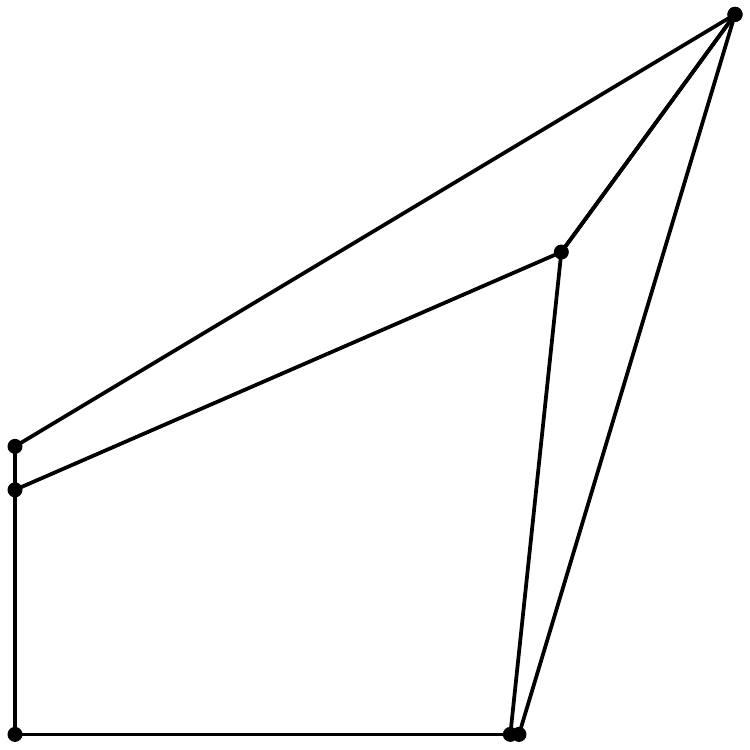}\includegraphics[scale=0.28, viewport=0cm -2cm 13cm 13cm, clip]{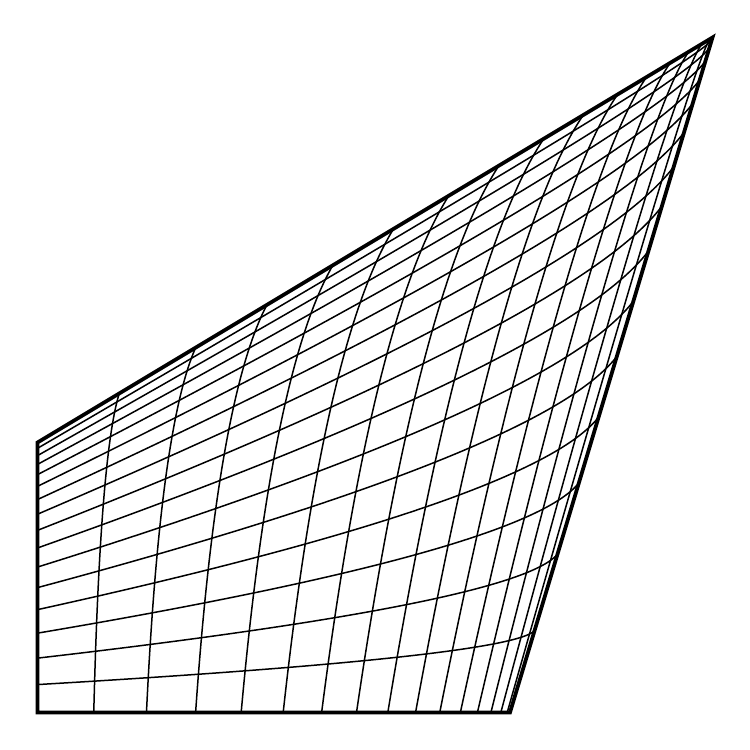}\label{Po_SingleQuadDeformCP}}
	\subfloat[$L^2$- and $H^1$-errors]{
		\includegraphics{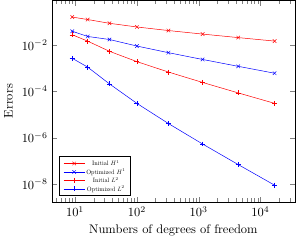}\label{Po_SingleQuadErrors}
	} 
	\caption{Poisson problem on a non-square quadrilateral domain with a corner singularity.}
\end{figure}

For the Poisson problem, we use the exact solution $u(x,y)=((x - 1)^2 + (y - 1)^ 2)^{\frac14}$ with a singularity at the upper right corner of the domain, see Figure~\ref{Po_SingleQuadSol}, and Dirichlet boundary conditions.
The deformed control net and the parameter lines are presented in Figure~\ref{Po_SingleQuadDeformCP}.
 We report the $H^1$- and $L^2$-errors in Figure~\ref{Po_SingleQuadErrors}.
 An improvement of almost four orders of magnitude is observed in the $L^2$-error while the $H^1$-error improves in more than one order of magnitude.

 In all examples up to now, we observed that the $L^2$-errors on the single triangles did not improve significantly after 10 iterations. Thus in the following we always perform 10 iterations without computing any $L^2$-errors on the individual
 triangles.

\subsection{Comparison of different sampling strategies}
In Section~\ref{sec:samplingh1}, we discussed different options for sampling input points for our method when optimizing for solving the Poisson problem. These correspond to sampling from the graph surface
\begin{linenomath}\begin{equation*}
	(x, y, u_{\mathit{init}}(x,y)) = (G_{\mathit{init}}(s,t), \hat u_{\mathit{init}}(s,t))
\end{equation*}\end{linenomath}
as well as from different directional derivatives $\frac{\partial}{\partial \vec v}\hat u_{\mathit{init}}$.

In this example we investigate the influence of this choice. In particular, we sample from $\hat u_{\mathit{init}}$, from the partial derivatives $\nabla \hat u_{\mathit{init}}\cdot e_i$
as well as from the directional derivatives  $\nabla \hat u_{\mathit{init}}\cdot Re_i$  for $i=1,2$. Here, $e_i$ are the canonical basis vectors and $R$ is the counter-clockwise rotation by $45^\circ$.
 
\begin{figure}
	\centering
	\subfloat[Exact solution]{\includegraphics[scale=0.18, viewport=3cm 3cm 36cm 25cm, clip]{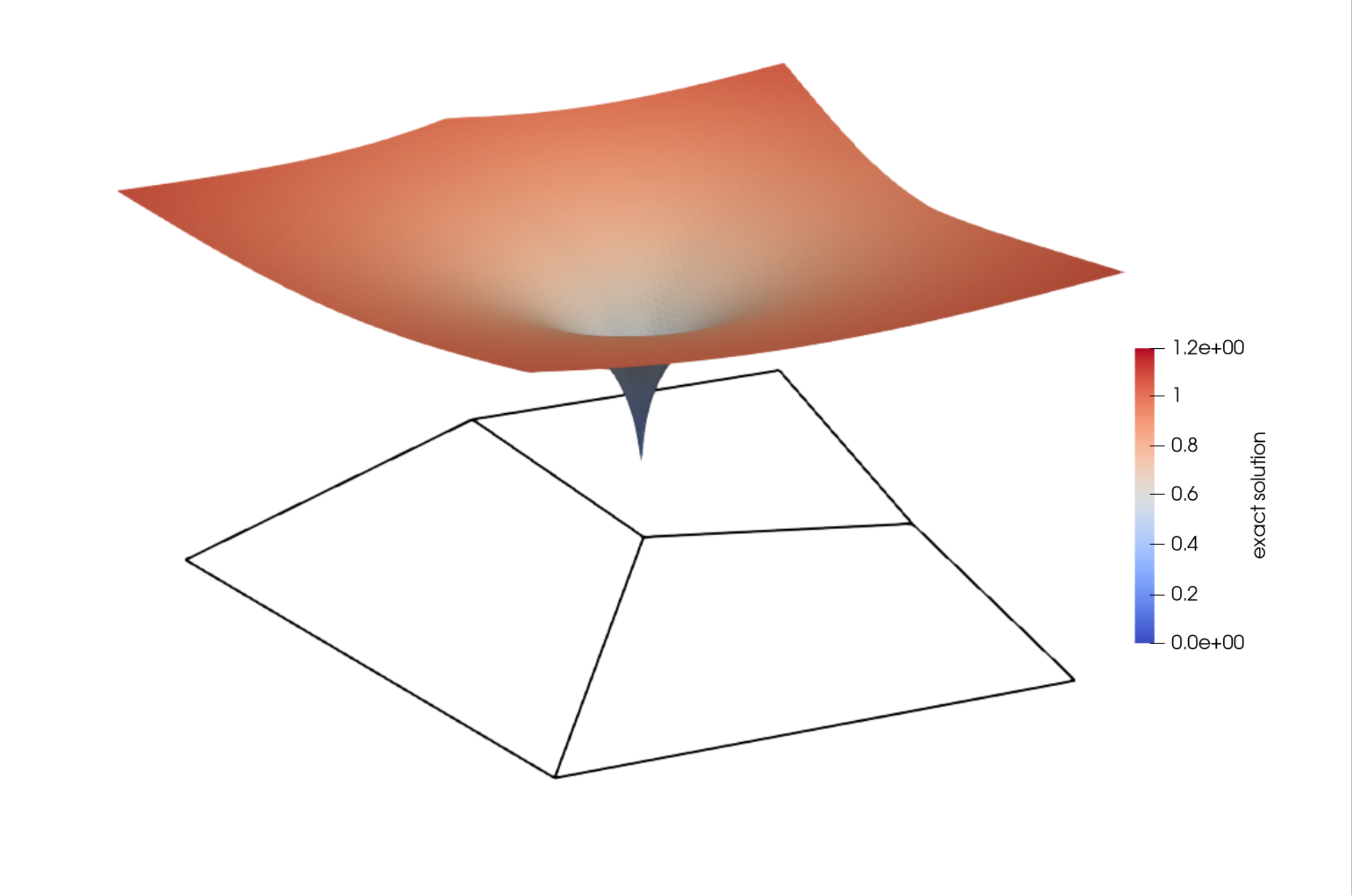}\label{Po_ThreePatchesSol}}
	\subfloat[$L^2$- and $H^1$-errors]{
		\includegraphics{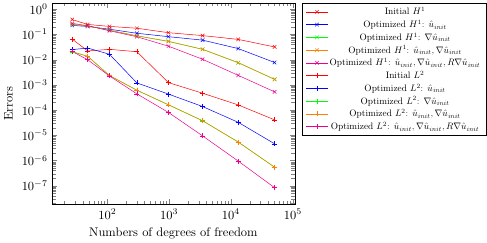}\label{Po_ThreePatchesErrors}}
\caption{Comparison of the different sampling strategies: Exact solution and convergence.}
\end{figure}

  To test the different choices, we solve the Poisson equation with Dirichlet boundary conditions on a multi-patch domain consisting of three patches with common vertex $(0, 0)$. The exact solution is chosen to be $u(x,y) = (x^2 + y^2 + 0.00001)^{\frac18}$, see Figure~\ref{Po_ThreePatchesSol}. Again, we perturb a singular function to obtain a function that is $C^\infty$ with very large $H^1$-norm.
  
\begin{figure}
	\centering
	\subfloat[Result for sampling from $\hat u_{\mathit{init}}$]{\includegraphics[scale=0.25, viewport=0cm 0cm 13.5cm 11.5cm, clip]{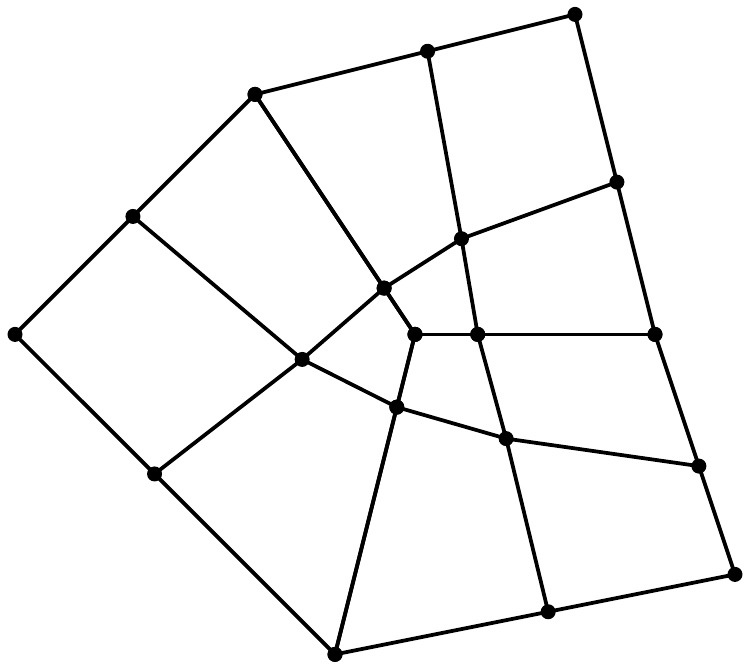}\includegraphics[scale=0.25, viewport=0cm 0cm 13cm 11cm, clip]{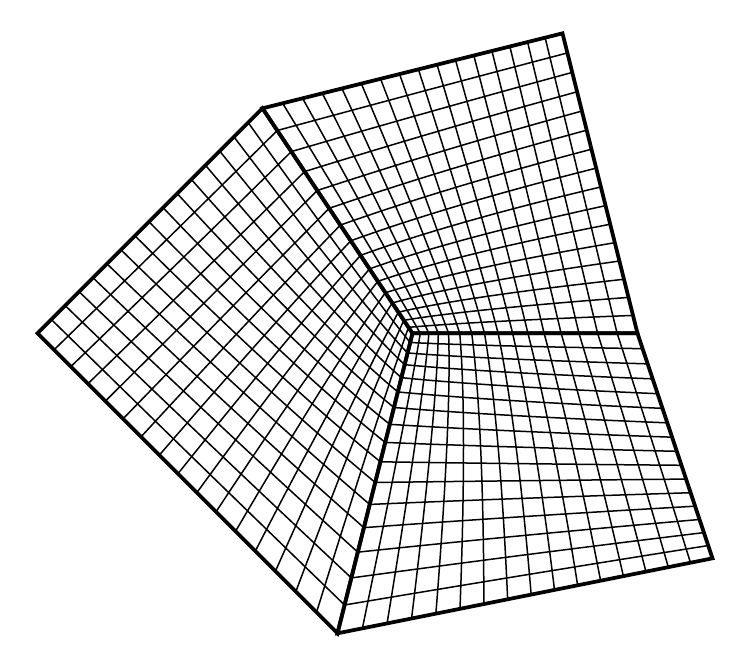}} \hspace{1cm}
	\subfloat[Result for sampling from $\nabla \hat u_{\mathit{init}}$]{\includegraphics[scale=0.25, viewport=0cm 0cm 13cm 11.5cm, clip]{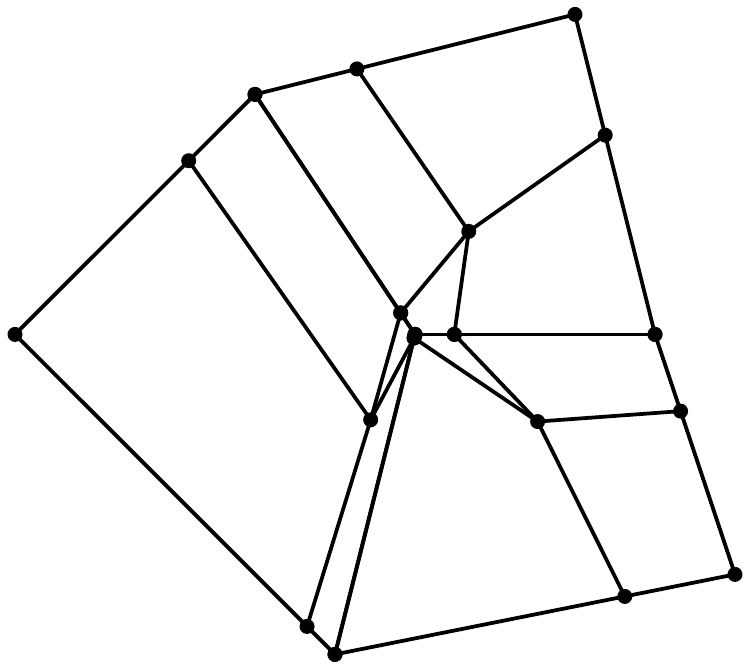}\includegraphics[scale=0.25, viewport=0cm 0cm 15cm 11cm, clip]{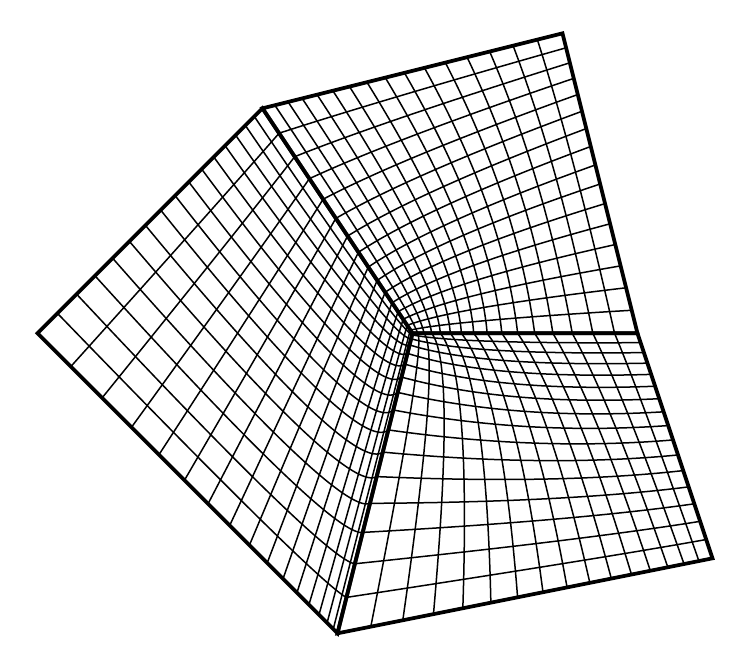}}\\
	\subfloat[Result for sampling from $\hat u_{\mathit{init}}$ and $\nabla \hat u_{\mathit{init}}$]{\includegraphics[scale=0.25, viewport=0cm 0cm 13.5cm 11.5cm, clip]{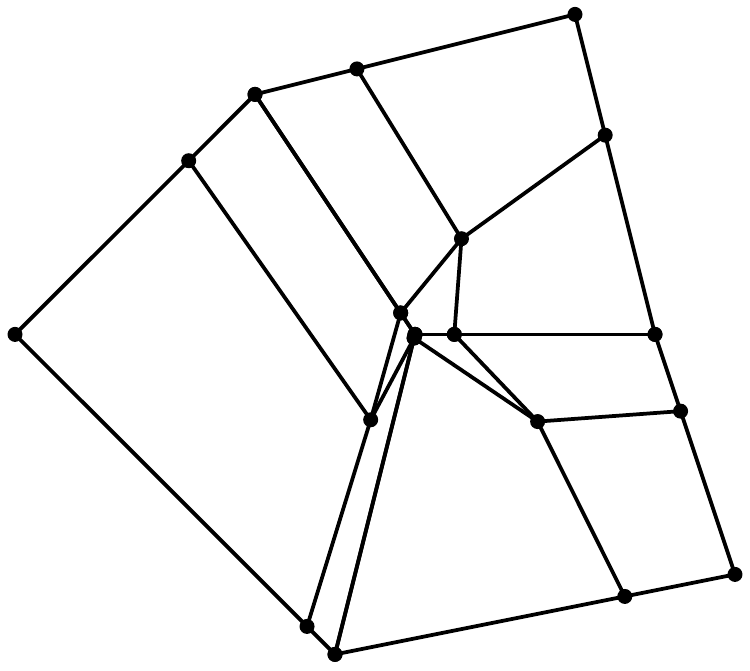}\includegraphics[scale=0.25, viewport=0cm 0cm 13cm 11cm, clip]{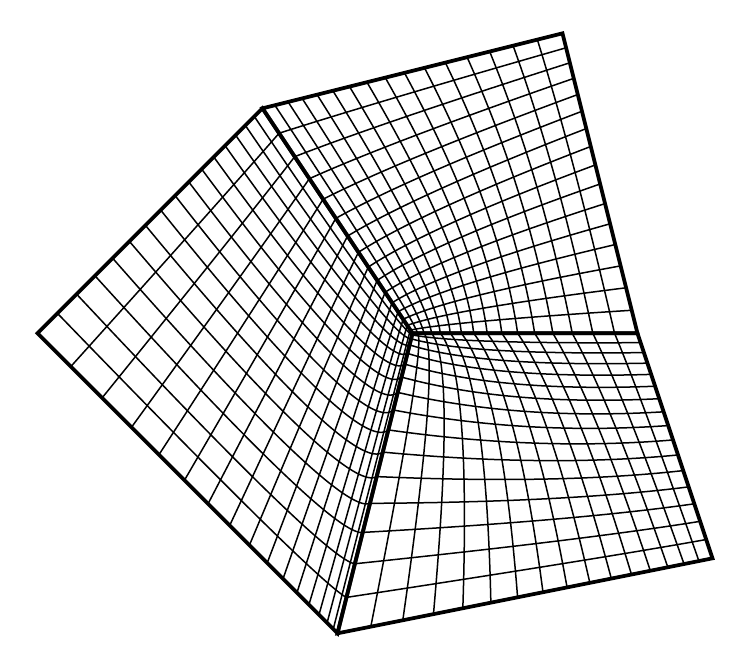}}
	\hspace{1cm}
	\subfloat[Result for sampling from $\hat u_{\mathit{init}}$, $\nabla \hat u_{\mathit{init}}$ and $R\nabla \hat u_{\mathit{init}}$]{\includegraphics[scale=0.25, viewport=0cm 0cm 13cm 11.5cm, clip]{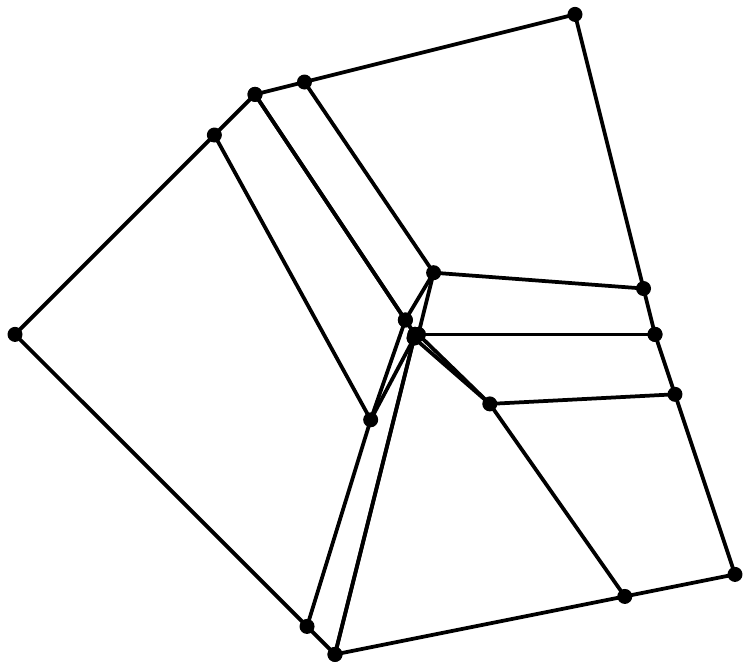}\includegraphics[scale=0.25, viewport=0cm 0cm 15cm 11cm, clip]{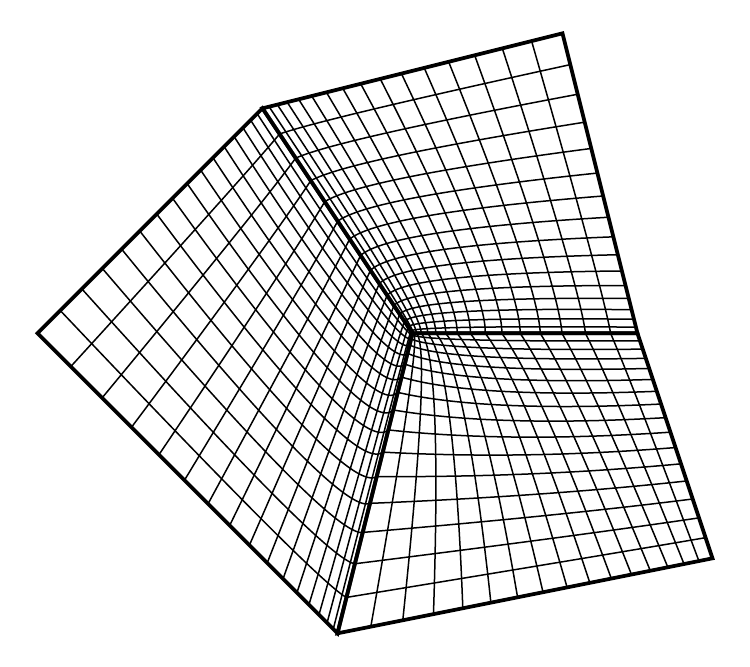}}\\
	\caption{Comparison of the different sampling strategies: Control nets and parameter lines for the different sampling strategies}\label{Po_ThreePatchesDeformCP}
\end{figure}

  In Figure~\ref{Po_ThreePatchesDeformCP} we show control nets and parameter lines for different sets containing these functions, always using the parameterization averaging. For each set of functions we choose the control points with maximum deformation, see Section~\ref{sec:samplingh1}.
 In Figure~\ref{Po_ThreePatchesErrors} we report the corresponding $H^1$- and $L^2$-error.
  One can observe that using the whole set of input functions leads to the largest improvement in the approximation error, both in the $L^2$-norm (more than two orders of magnitude) and in the $H^1$-seminorm (almost two orders of magnitude). The improvement for the smaller sets of functions is smaller, but still significant.
  
  We conclude that it is beneficial to use all five functions for sampling and therefore report only the result of this strategy for all other examples where the Poisson equation is solved, including the previous examples. Note that for $L^2$-projection only the graph surface of $\hat u_{\mathit{init}}$ is used for sampling.

\subsection{Solutions with multiple singularities}
\begin{figure}
	\centering
	\subfloat[Exact solution]{\includegraphics[scale=0.18, viewport=11cm 0cm 32cm 22cm, clip]{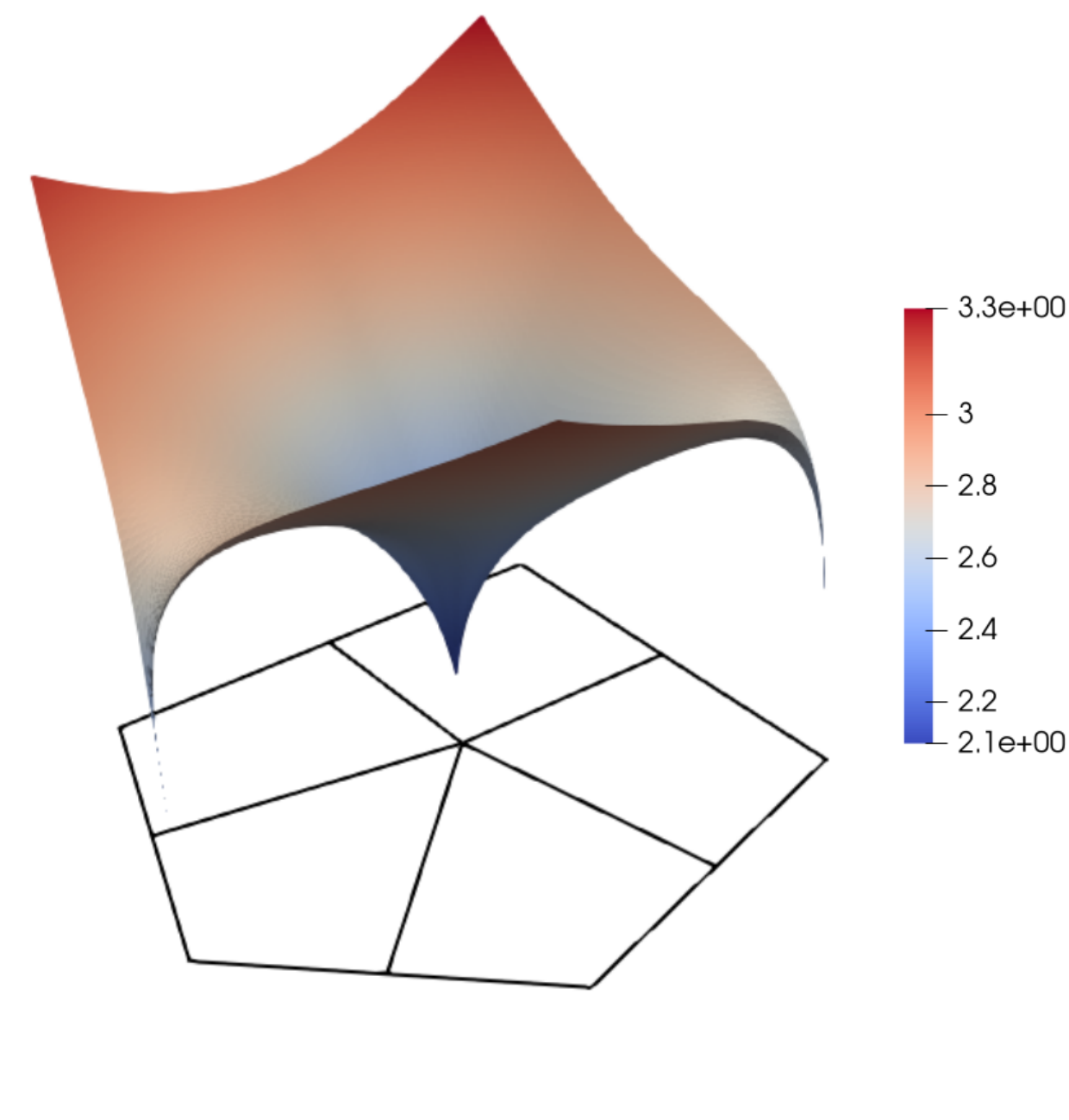}\label{Po_PentagonSol}}
	\subfloat[Resulting control net and parameter lines]{\includegraphics[scale=0.28, viewport=-1cm -2cm 13.5cm 13cm, clip]{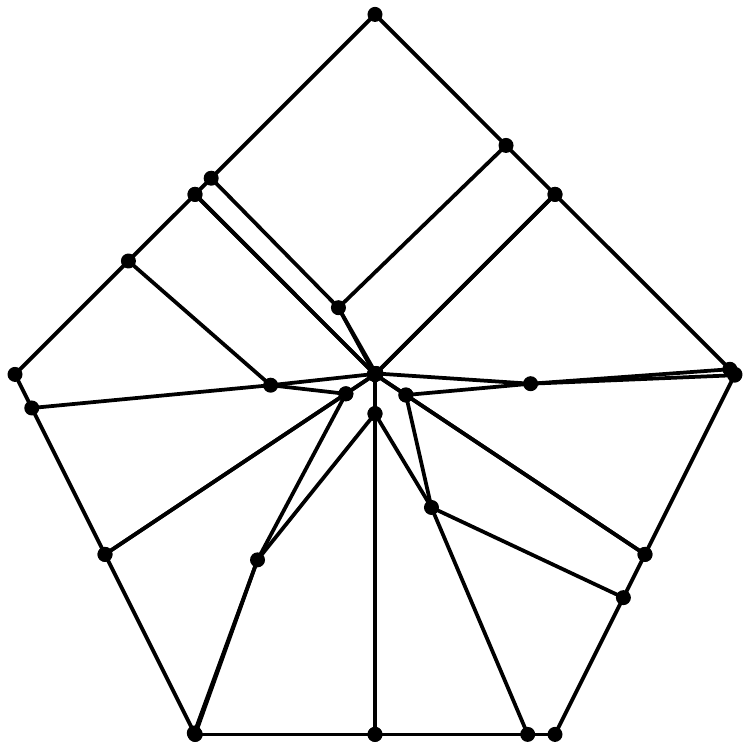}\includegraphics[scale=0.28, viewport=0.5cm -2cm 13cm 15cm, clip]{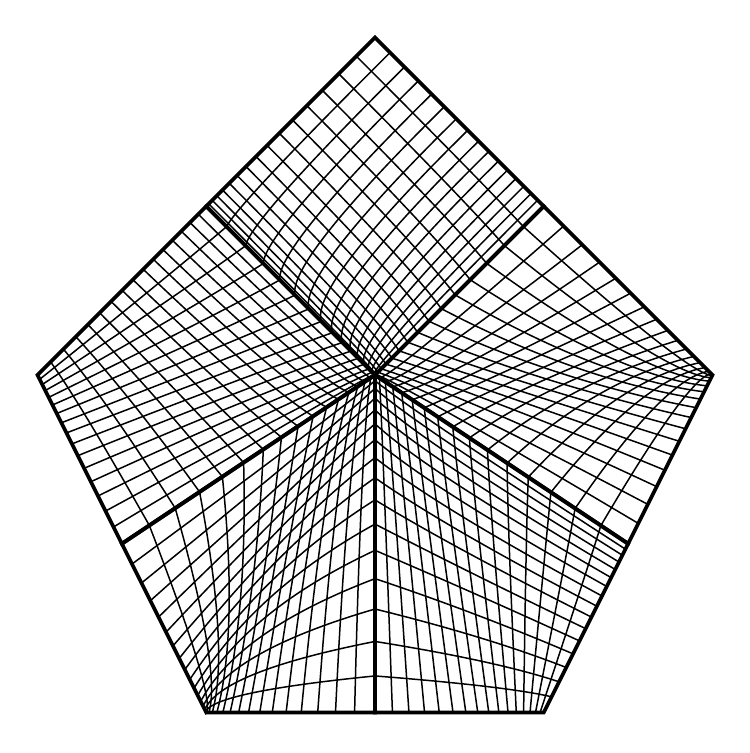}\label{Po_PentagonDeformCP}}
	\subfloat[$L^2$- and $H^1$-errors]{
		\includegraphics{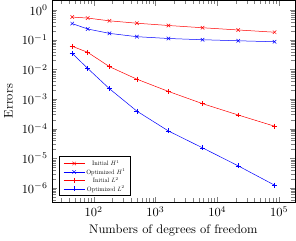}\label{Po_PentagonErrors}}
	\caption{Poisson problem with a solution with multiple singularities: Exact solution, deformed control net and parameter lines, convergence plot.}
\end{figure}

In this section we consider two Poisson problems with Dirichlet boundary conditions with exact solutions having multiple point singularities and two sides singularity. In each example we use the parameterization averaging to obtain the deformed control net. 
\paragraph{Multiple point singularities on a pentagon}
We start with the exact solution $u(x,y)=((x + 0.5)^2 + (y + 1)^2)^{\frac18} + (x^2 + y^2 + 0.00001)^{\frac14} + ((x - 1)^2 + y^2)^{\frac18}$ on a pentagon that is represented by five bilinear quadrilateral patches with a common vertex at $(0,0)$, see Figure~\ref{Po_PentagonSol}.
At $(0,0)$, this function is close to singular, the two singularities at $(-\frac{1}2,-1)$ and $(1,0)$ are at corners of two patches of the domain.

\begin{figure}
	\centering
	\subfloat[Exact solution]{\includegraphics[scale=0.18, viewport=12cm 1cm 28cm 25cm, clip]{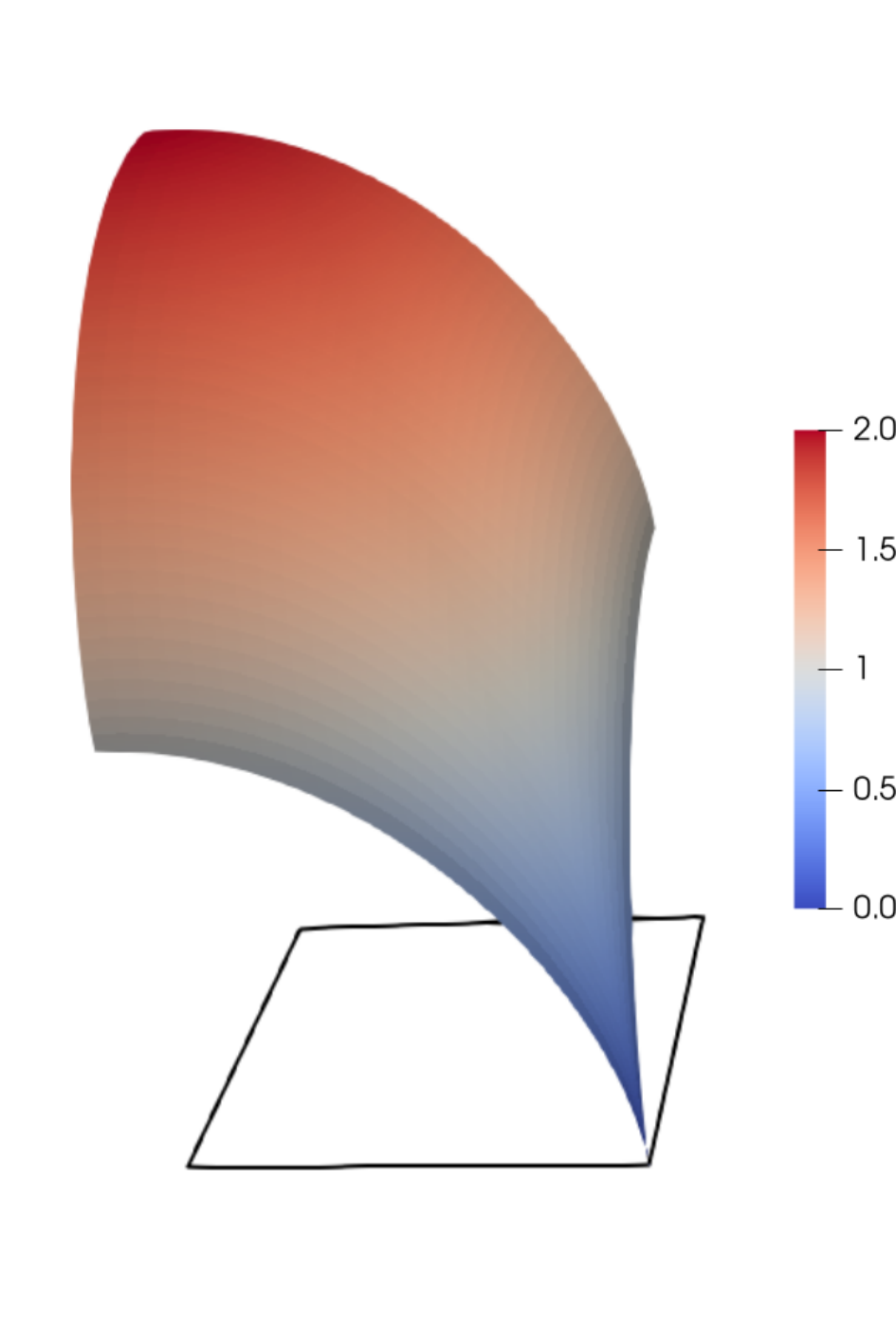}\label{Po_SingleSquareTwoSidesSingSol}}
	\subfloat[Resulting control net and parameter lines]{\includegraphics[scale=0.27, viewport=-1cm -2cm 13.5cm 13cm, clip]{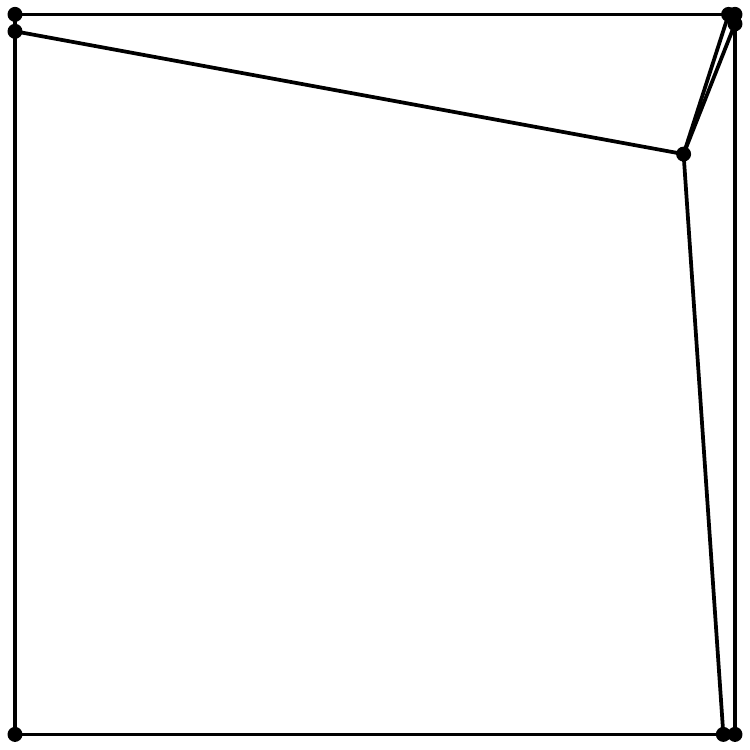}\includegraphics[scale=0.28, viewport=0cm -2cm 13cm 15cm, clip]{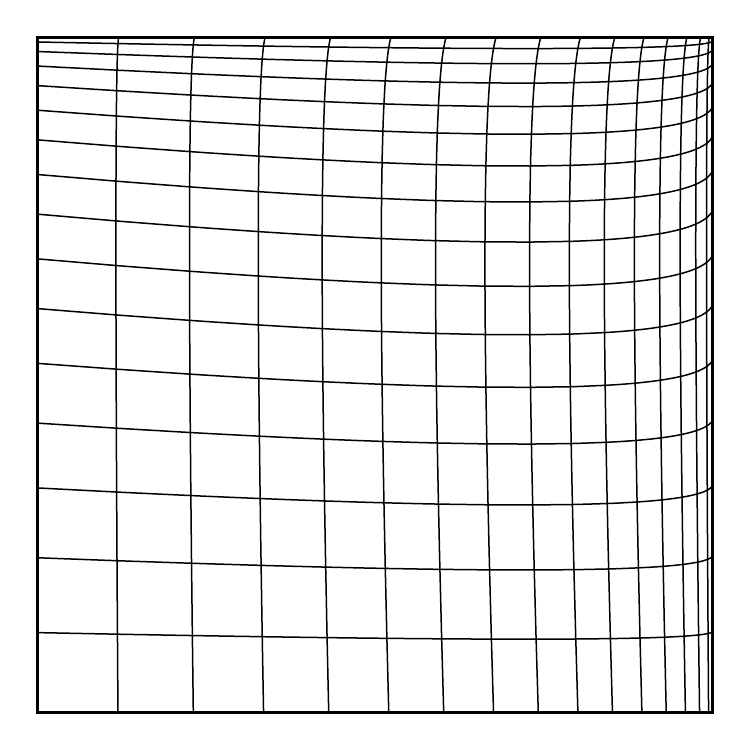}\label{Po_SingleSquareTwoSidesSingDeformCP}}
	\subfloat[$L^2$- and $H^1$-errors]{
		\includegraphics{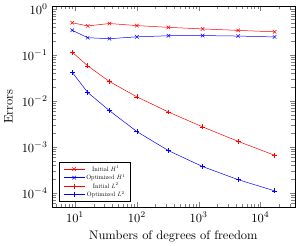}\label{Po_SingleSquareTwoSidesSingErrors}}
	\caption{Poisson problem with a solution with multiple side singularities: Exact solution, deformed control net and parameter lines, convergence plot.}
\end{figure}

The deformed control net and the parameter lines after optimizing the parameterization using our method are shown in Figure~\ref{Po_PentagonDeformCP}. We observe that the control points move towards all three singularities.
We report the $L^2$- and $H^1$-errors in Figure~\ref{Po_PentagonErrors}.
 An improvement of more than two orders of magnitude is observed in the $L^2$-error. The improvement of the $H^1$-error is smaller but after only two steps of refinement we achieve a smaller $H^1$-error than after seven refinement steps of the original parameterization.
\paragraph{Multiple side singularities}
Now we consider an exact solution with two side singularities, $u(x,y)=(1 - x^2)^{\frac35} + (1 - y^2)^{\frac35}$, on the unit square, see Figure~\ref{Po_SingleSquareTwoSidesSingSol}. Such an example is in the spirit of flow problems with boundary layers, which can appear on several edges simultaneously, cf.~\cite{hughes2005isogeometric}. This function is singular along the top and the right edge of the square. Using classical methods, adapting the discretization correctly near the corner between these edges poses a significant difficulty.
The deformed control net with the corresponding parameter lines obtained from the parameterization averaging are presented in Figure~\ref{Po_SingleSquareTwoSidesSingDeformCP} and we report the $L^2$- and $H^1$-errors after several steps of $h$-refinement in Figure~\ref{Po_SingleSquareTwoSidesSingErrors}. One can observe an improvement of almost one order of magnitude of the $L^2$-error as well as an improvement of the $H^1$-error, especially for coarse meshes. 

\subsection{Multi-patch isogeometric analysis on domains with a reentrant corner}
When considering a Poisson problem on a domain with a single reentrant corner, the solution $u$ is in general of low regularity. We have that $u\in H^1(\Omega)$ is the sum of a regular part $u_r \in H^\ell(\Omega)$, with $\ell\geq 2$, and a singular part $u_s \in H^{1+\epsilon}(\Omega)$, for some $\epsilon >0$. In case of a reentrant corner with angle $\omega \in \left]\pi,2\pi\right[$, we actually know that the singular part satisfies $u_s = \xi(r) r^\lambda \sin(\lambda \theta)$, when parameterized in polar coordinates, with $0 < r < R$, $0 < \theta < \omega$, $\lambda = \pi/\omega$ and $\xi(r)$ being a $C^\infty$-function that has support only in the region $0<r<R$. We thus have $u_s \in W^{2,\frac{2}{2-\lambda}}(\Omega) \subseteq H^{1+\lambda}(\Omega)$. In case of an L-shape, with $\omega = 3\pi/2$, we thus have $u_s \in H^{5/3}$.

While away from the singular point a local error estimate of the form
\[
 | u - \Pi_h u |_{H^1(Q)} \lesssim h_Q^{2} \| u \|_{H^{3}(\tilde Q)}
\]
holds for $p=2$, where $Q$ is an element in physical space, $h_Q$ its size and $\tilde Q$ its support extension, cf.~\cite{bazilevs2006isogeometric}, no such estimate can be shown for elements near the singularity. Thus, the mesh must be graded, such that the local error from the singular function is balanced with the local errors away from the singularity.

\begin{figure}
	\centering
	\subfloat[Exact solution]{\includegraphics[scale=0.16, viewport=10cm -1cm 31.6cm 24cm, clip]{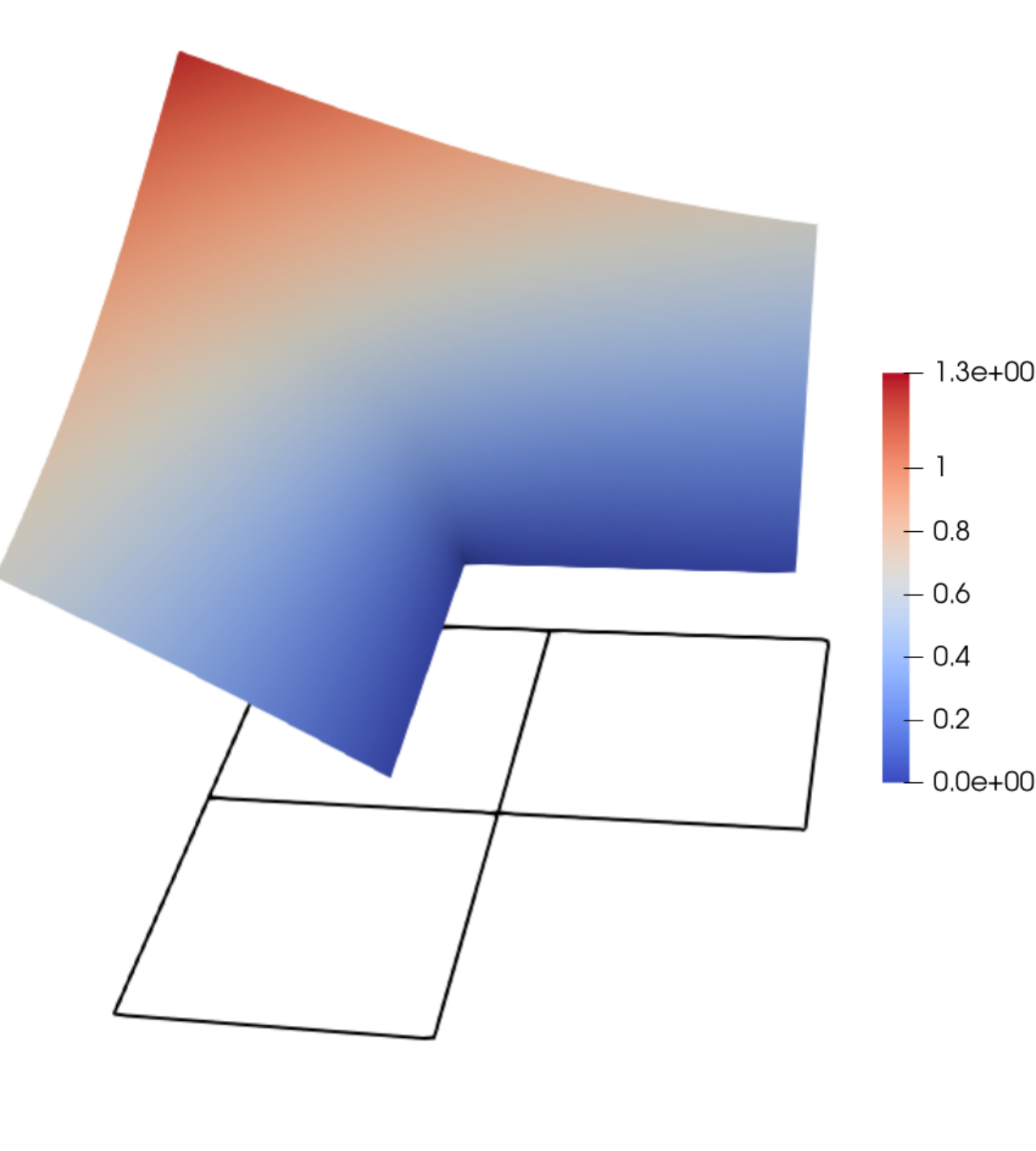}\label{Po_LshapeStatHeatSol}}
	\subfloat[Resulting control net and parameter lines]{\includegraphics[scale=0.27, viewport=-0.8cm -2cm 13.5cm 13cm, clip]{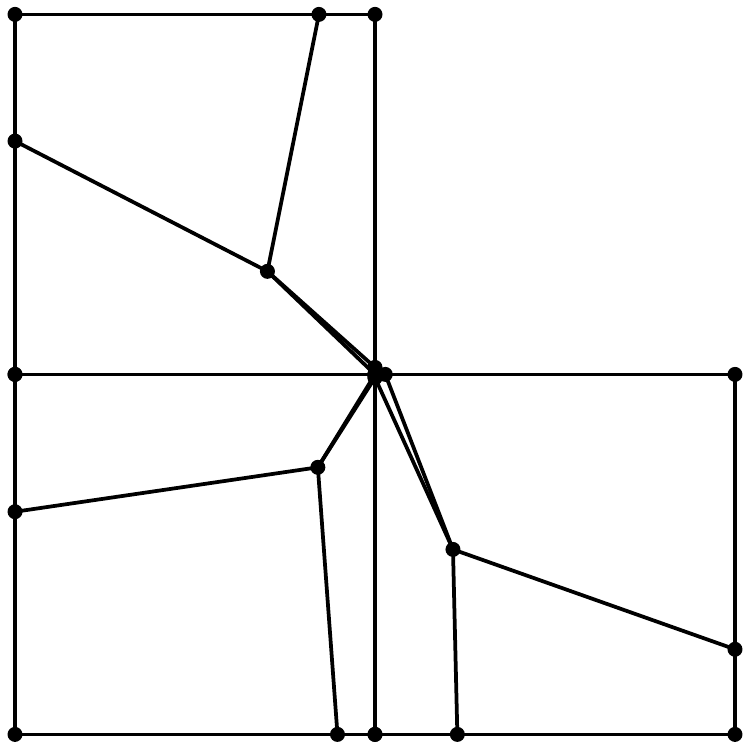}\includegraphics[scale=0.3, viewport=0.6cm -1.5cm 12.5cm 15cm, clip]{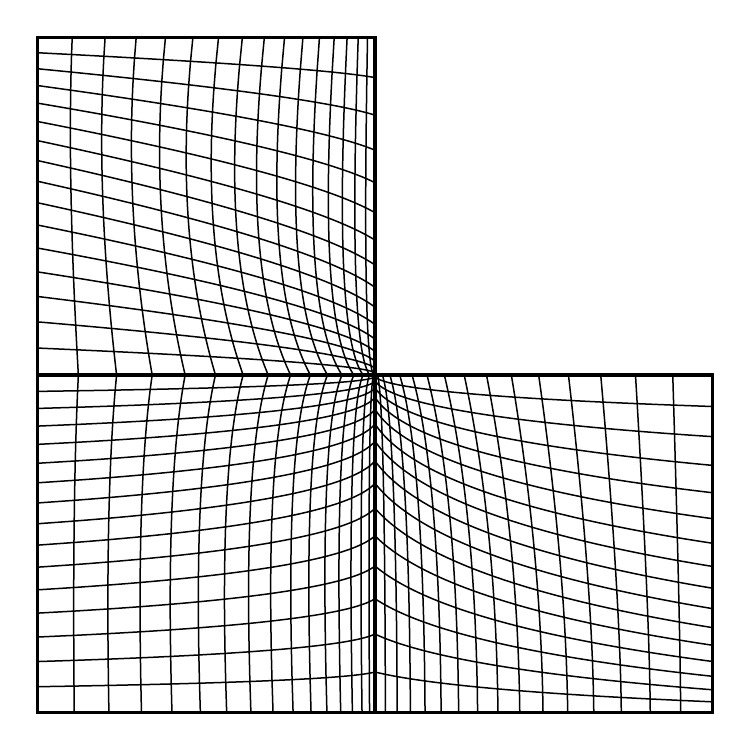}
		\label{Po_LshapeStatHeatDeformCP}} 
	\subfloat[$L^2$- and $H^1$-errors]{
		\includegraphics{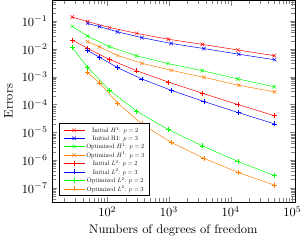}\label{Po_LshapeStatHeatErrors}}
	\caption{Stationary solution of the heat equation on an L-shape domain: Exact solution, deformed control net and parameter lines, convergence plot.}\label{fig:L-shape-heat-all}
\end{figure}

\paragraph{Stationary heat conduction on an L-shape domain}
In an example  taken from~\cite{dorfel2010adaptive}, we compute the stationary solution to the heat equation on the L-shape domain $[-1,1]^2\setminus[0,1]^2$ with exact solution
\begin{linenomath}\begin{equation*}
u(r,\phi) = r^{\frac23}\sin(\frac{2\phi-\pi}{3})
\end{equation*}\end{linenomath}
in polar coordinates. This functions satisfies $\Delta u = 0$ in Cartesian coordinates. We impose zero Dirichlet boundary conditions on the two boundary segments emerging from the reentrant corner, see Figure~\ref{Po_LshapeStatHeatSol}. On the remaining boundary we impose Neumann boundary conditions obtained from the exact solution.

In Figure~\ref{Po_LshapeStatHeatDeformCP}, we show the resulting control points and parameter lines. 
We solved the equation both in the isogeometric function space generated by biquadratic tensor-product splines and in the space generated by bicubic tensor-product splines.
For both degrees, our method results in a large improvement of the errors, both in the $L^2$- and the $H^1$-norms, see Figure~\ref{Po_LshapeStatHeatErrors}. As expected, the bicubic splines perform only slightly better than the biquadratic splines, due to the reduced Sobolev regularity of the solution.
Comparing our results for degree 3 with the results reported in~\cite{dorfel2010adaptive}, where the authors solved the same problem by adaptive refinement using cubic T-splines, we note that we achieve slightly smaller errors in the $L^2$-norm using similar numbers of degrees of freedom. The errors in the $H^1$-norm seem to be close for similar numbers of degrees of freedom.
Note that our method maintains the tensor-product structure, which is beneficial in terms of efficiency and enables us to compute many steps of refinement.

\begin{figure}
	\centering
	\subfloat[Discrete solution]{\includegraphics[scale=0.15, viewport=12cm -1cm 30.5cm 24cm, clip]{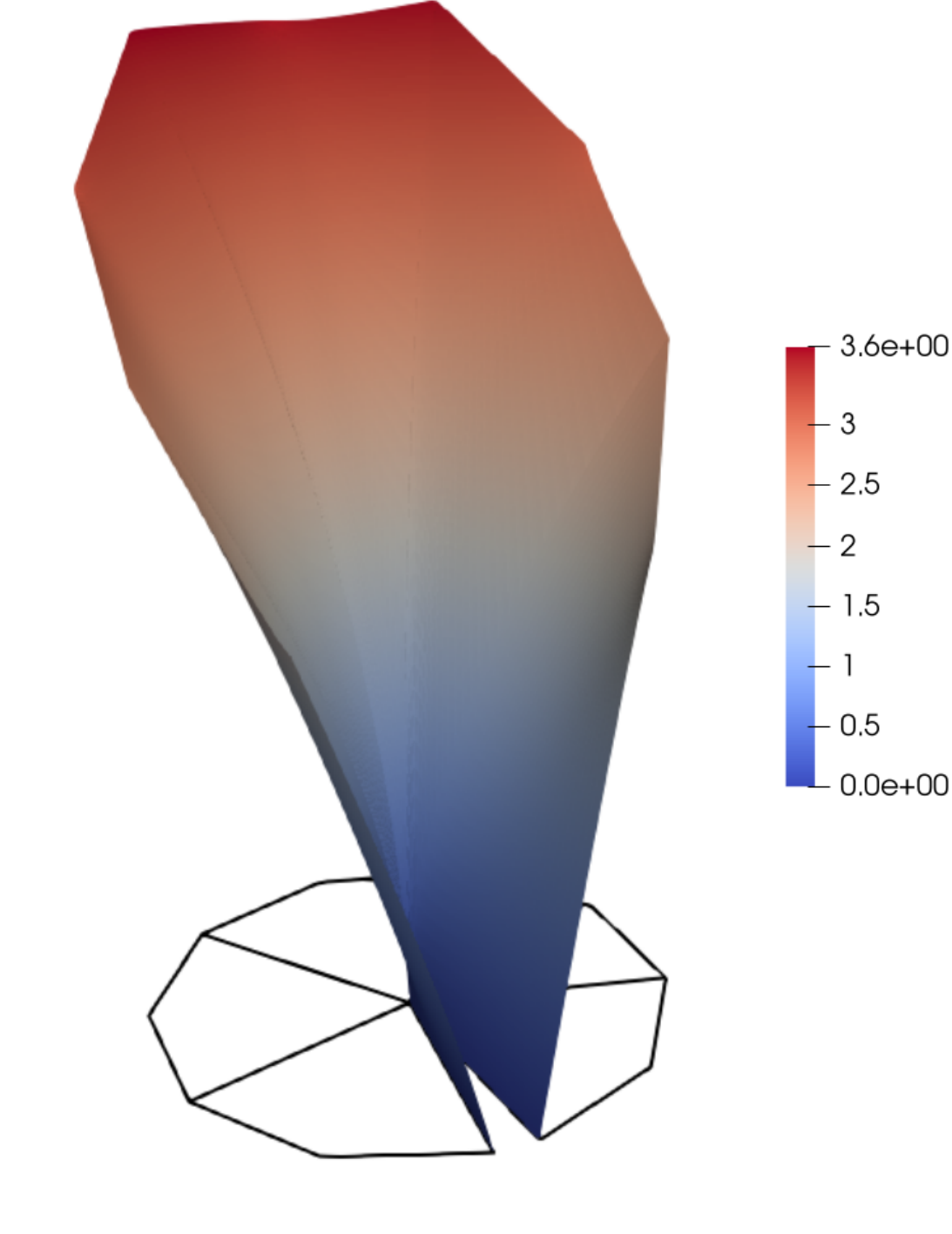}\label{Po_HexagonNeumSol}}
	\subfloat[Resulting control net and parameter lines]{\includegraphics[scale=0.3, viewport=-1cm -2cm 13.5cm 12.5cm, clip]{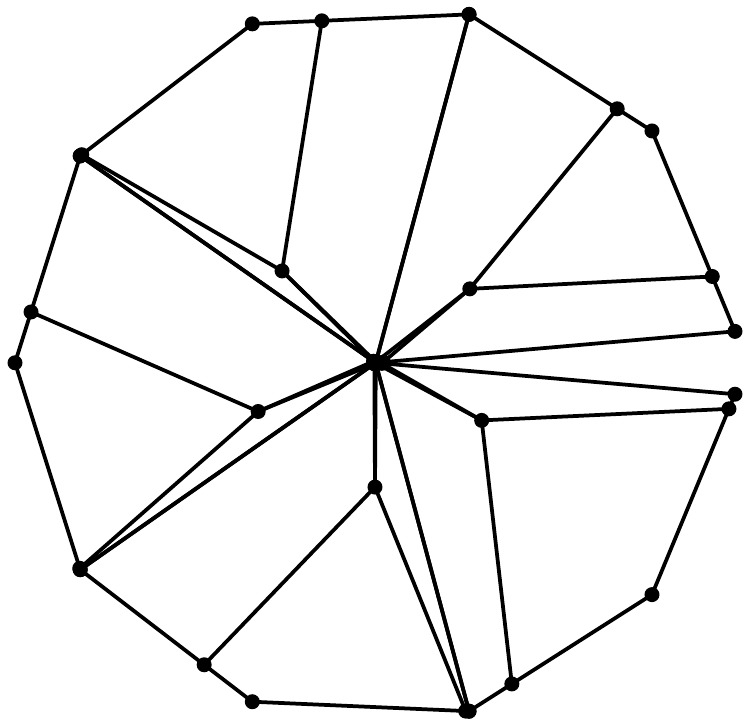}\includegraphics[scale=0.33, viewport=0.5cm -1.5cm 13cm 12cm, clip]{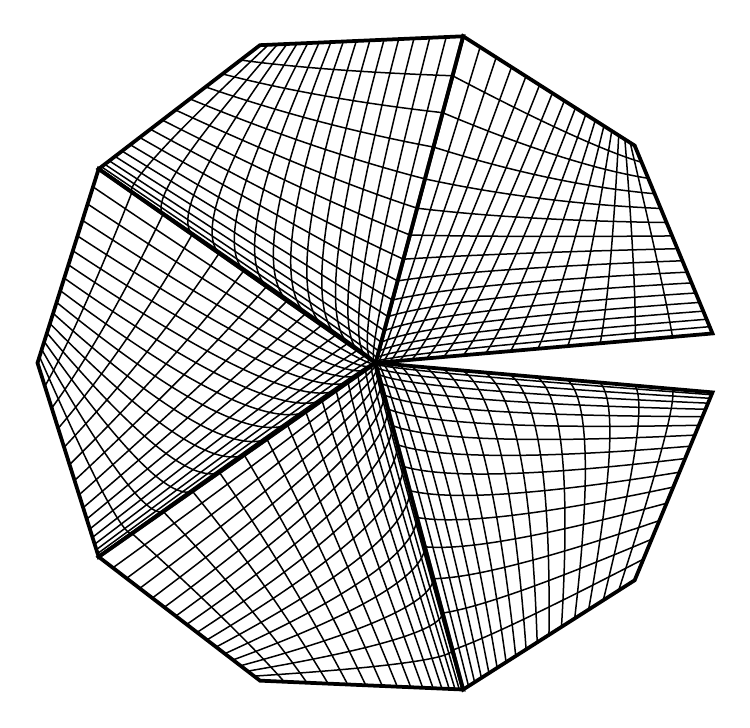}\label{Po_HexagonNeumDeformCP}}
	\subfloat[$\max_\Omega(\|\nabla u_h\|)$]{
		\includegraphics{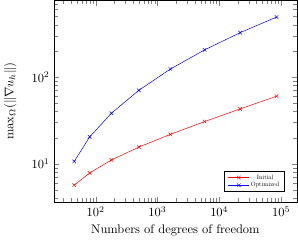}\label{Po_HexagonNeumErrors}}
	\caption{Poisson problem on a cracked disk: Discrete solution $u_h$ after 7 refinement steps, deformed control net and parameter lines, plot of the maximum of $\|\nabla u_h\|$ depending on $h$.}
\end{figure}

\paragraph{Crack singularity}
Finally, we consider a Poisson problem without known exact solution, inspired by a crack simulation. To this end, we construct a multi-patch domain consisting of five patches meeting at a common vertex with an reentrant corner with angle $\omega = 350\degree$. We impose zero Dirichlet conditions on the edges emerging from the reentrant corner and Neumann boundary conditions $g_N=1$ everywhere else. Figure~\ref{Po_HexagonNeumSol} shows a numerical approximation to the solution.
In Figure~\ref{Po_HexagonNeumDeformCP} we show the control points and parameter lines resulting from our method. As we see, the control points moved towards the singularity.
Since we do not have a reference solution and the approximate solution generated by uniform refinement of the initial parameterization corresponds to a different spline space than the one resulting from the deformed parameterization, we do not estimate the errors. Instead, we plot the maximum of the norm of the gradient of the approximate solution, obtained by sampling at the quadrature nodes. As expected, the approximate solutions diverge towards infinity, with our method diverging significantly faster.

\section{Possible extensions and conclusion}\label{sec:extension-conclusions}

In this paper we propose a strategy for $r$-adaptivity of isogeometric discretization spaces, where a domain reparameterization is found with the help of a neural network. The strategy requires that the domain is first segmented into (bilinear) patches, which are then reparameterized independently. Standard $h$-refinement over the reparameterized patches then results in suitably graded meshes. This approach has several advantages compared to standard, mesh refinement based, adaptivity. The resulting spaces are still patch-wise tensor-products, which allows for the use of very efficient quadrature and assembly routines. The application of a neural network is very efficient, saving time and computational effort that would otherwise be needed for error estimation and local refinement. The neural network requires only data from the discrete solution on a coarse mesh and no other information of the PDE problem.

The proposed approach is flexible and was tested on a large variety of model configurations. The obtained results show that the approach can be used as is, but improvements can be done nonetheless. In the future we plan on expanding the method in several directions. The initial parameterization should be allowed to be a general multi-patch spline mapping. This can be achieved e.g. by composing the patch parameterization with a reparameterization of the parameter domain onto itself. Currently, we restrict ourselves to biquadratic reparameterizations. There is no inherent difficulty to extend the approach to general degrees. However, this requires training a new network. Moreover, the reparameterization should allow for changes of the patch interfaces, while keeping only the global domain boundary fixed. Such a change would require a different interpretation of the patch segmentation and a different setup for computing the reparameterization. As a next step, we intend to improve upon the strategy by alternating between local reparameterizations and global $h$-refinement. Such a strategy would generate a (possibly low rank) reparameterization that maintains the tensor-product structure and improves the approximation even further.

\section*{Acknowledgments}
The authors were supported by the Linz Institute of Technology (LIT) and the government of Upper Austria through the project LIT-2019-8-SEE-116 entitled ``PARTITION -- PDE-aware isogeometric discretization based on neural networks''. This support is gratefully acknowledged.

\bibliography{refs}

\end{document}